\documentclass[12pt,reqno,nosumlimits]{amsart}
\usepackage{amsfonts}

\input{xyv2.tex}
\usepackage{color}

\newtheorem{theorem}{Theorem}[section]
\newtheorem{proposition}[theorem]{Proposition}
\newtheorem{lemma}[theorem]{Lemma}
\newtheorem{corollary}[theorem]{Corollary}
\theoremstyle{definition}
\newtheorem{definition}[theorem]{Definition}
\newtheorem{example}[theorem]{Example}

\theoremstyle{remark}
\newtheorem{remark}[theorem]{Remark}
\newtheorem{noname}[theorem]{\bf\hspace*{-1ex}}

\def\Id{\mathrm{Id}}
\def\Pd{P_\ast}
\def\ama{{}_A\mathfrak{M}_A}
\def\bmb{{}_B\mathfrak{M}_B}
\def\calE{\mathcal{E}}
\def\calM{\mathfrak{M}}
\def\cmc{{}_C\mathfrak{M}_C}
\def\cmm{\mathcal{CM}_m(H)}
\def\cm{\mathcal{CM}(H)}
\def\dd{d_\ast}
\def\gr{\mathrm{gr}}

\def\dwotB{\widehat{\underset{B}{\otimes}\,}}
\def\otB{\otimes_B}

\newcommand{\adjoint}{-\hspace{-1.2ex}-\hspace{-0.75ex}-\hspace{-1.1ex}\mid}
\newcommand{\Coker}{\mathrm{Coker}}

\newcommand{\Hom}{\mathrm{Hom}}
\newcommand{\Ind}{\mathrm{Ind}}

\newcommand{\Tor}[2]{\mathrm{Tor}^{#1}_{#2}}
\newcommand{\Tot}{\mathrm{Tot}}

\newcommand{\an}[1]{_{\langle{#1}\rangle}}
\newcommand{\ccd}{\mathrm{CC}_{\ast\ast}}

\newcommand{\de}[1]{_{({#1})}}

\newcommand{\mm}{\mathfrak{M}}
\newcommand{\pr}[3]{{\textstyle{\prod}}_{j=#1}^{#2}a^j_{\langle
#3\rangle}}
\newcommand{\rmd}[2]{\mathrm{#1}_{#2}}
\newcommand{\rmu}[2]{\mathrm{#1}^{#2}}
\newcommand{\tp}[1]{\otimes_{#1}}
\newcommand{\tsum}{{\textstyle\sum}}
\newcommand{\wccd}{\widetilde{\mathrm{CC}}_{\ast\ast}}
\newcommand{\wotB}{\widehat{\otimes}_B}
\newcommand{\wrmd}[2]{\widetilde{\mathrm{#1}}_{#2}}
\newcommand{\yd}{\mathcal{YD}(H)}
\newcommand{\Ulie}[1]{{U}(\mathfrak{#1})}
\newcommand{\Uf}[1]{{U}_f(\mathfrak{#1})}
\newcommand{\Lie}[1]{\mathfrak{#1}}
\def\vsu{\vspace*{1ex}}\def\vsd{\vspace*{2ex}}\def\vst{\vspace*{3ex}}
\begin{document}
\title{Cyclic homology of Hopf Galois extensions and Hopf algebras}

\begin{abstract}
Let $H$ be a Hopf algebra and let $\mathcal{CM}_m(H)$ be the category of all left
$H$--modules and right $H$--comodules satisfying the following two compatibility
relations:
\begin{eqnarray*}
    &&\rho(hm)=\tsum h\de{2}m\an{0}\otimes h\de{3}m\an{1}Sh\de{1},\quad \text{for all }m\in
        M \text{ and }h\in H.\\
    &&\tsum m\an{1}m\an{0}=m,\quad \text{for all }m\in M.
\end{eqnarray*}
An object in $\mathcal{CM}_m(H)$ will be called a modular crossed
module (over $H$). For example, if $A/B$ is an $H$--Galois
extension, then the quotient $A_B:=A/[A,B]$ of $A$ modulo
commutators $[A,B]$ is a modular crossed module over $H$. More
particularly, $H$ itself can be regarded as an object
$_{ad}H\in\mathcal{CM}_m(H)$.

The category $\mathcal{CM}_m(H)$ has very nice homological
properties: it is abelian and contains enough injective objects.
Furthermore, if $K$ is a Hopf subalgebra of $H$, then   the
categories $\mathcal{CM}_m(K)$ and $\mathcal{CM}_m(H)$ can be
related by a functor $\mathrm{Ind}_K^H(-)$, where
$\mathrm{Ind}_K^HM:=H\otimes_KM$, with appropriate structures.

 To every $M\in\mathcal{CM}_m(H)$ we associate a cyclic object
$\mathrm{Z}_\ast(H,M)$. The cyclic homology of
$\mathrm{Z}_\ast(H,M)$ extends the usual cyclic homology of the
algebra structure of $H$ (for $M:={}_{ad}H$). The relative cyclic
homology of an $H$--Galois  extension $A/B$ can be also regarded
as a particular case ($M:=A_B$).

We compute the cyclic homology of $\mathrm{Ind}_K^HM:=H\otimes_KM$
when $K$ is cocommutative, and $M$ decomposes as a direct sum of
one dimensional subcomodules such that the associated group--like
elements are central. As a direct application of this result, we
describe the relative cyclic homology of strongly graded algebras.
In particular, we calculate the (usual) cyclic homology of group
algebras and quantum tori.

Finally, when $H:=\Ulie g$ is the enveloping algebra of a Lie
algebra $\Lie g$, we construct a spectral sequence that converges
to the cyclic homology of $H$ with coefficients in an arbitrary
modular crossed module $M$. We also show that the cyclic homology
of almost symmetric algebras is isomorphic to the cyclic homology
of $H$ with coefficients in a certain modular crossed--module.
\end{abstract}

\author{}

\author{Pascual Jara$^{\dag}$}
\thanks{$^{\dag}$ Research partially supported by DGES BMF2001--2823
and
       FQM-266 (Junta de Andaluc{\'\i}a Research Group).}
\address{Department of Algebra,
         University of Granada,
        18071--Granada, Spain}
\email{pjara@ugr.es}%

\author{Drago\c s \c{S}tefan$^\ddag$}
\thanks{$^\ddag$ Research supported by Grant SB2000-0264 from the Spanish
Ministry of Education.}
\address{University of Bucharest,
         Faculty of Mathematics,
         Str. Academiei 14, RO-70109,
         \hspace*{2.75ex}Bucharest, Romania}
\email{dstefan@al.math.unibuc.ro}%
\date{\today}

 \maketitle
 \pagestyle{myheadings}

\section*{Introduction}
\markboth{\sc Introduction}{\sc Introduction}

Cyclic cohomology was invented by A. Connes as a replacement of the de Rham cohomology
of varieties, \cite{Co1}. Since then, a lot of work has been done in order to compute
the cyclic (co)homology of certain classes of algebras. We would like to recall here
only some of the results obtained, namely those that somehow are connected to our
paper.

The cyclic homology of group algebras over fields of
characteristic $0$ was computed by Burghelea, \cite{Bu}. For a
complete algebraic proof of Burghelea's result the reader is
referred to \cite{Mar}, while a relative variant of this
computation can be found in \cite{Scha}.

Crossed products (with trivial cocycle) are generalizations of
group algebras. Let $B$ be an algebra on which a group $G$ acts by
algebra automorphisms. Let $A$ be the free left $B$ module having
a basis $\{e_g\mid g\in G\}$. One can define an algebra structure
on $A$ by setting
\[
(xe_g)\cdot(ye_h):=x{}(g.y)e_{gh},
\]
where $x,y\in B$, and $g.y$ denotes the action of $g$ on $y$.
Cyclic homology of crossed products was calculated by Feigin and
Tsygan, \cite{FT}. For related work on this case see also
\cite{GJ,Ni}.

The cyclic homology of enveloping algebras of Lie algebras is also
known. As a matter of fact, C. Kassel in \cite{Ka} described the
cyclic homology of all almost symmetric algebras ($A$ is called
almost symmetric if it is non--negatively filtered such that
$\gr{A}$ is the symmetric algebra of $\gr_1{A}$). The computation
of cyclic homology of $\Ulie{g}$ was also performed in \cite{FT}.

Cyclic (co)homology of Hopf algebras was introduced by Connes and
Moscovici in order to compute the index of transversally elliptic
operators of foliations. To every Hopf algebra $H$ and every
modular pair in involution $(\sigma, \delta)$ they associated a
cocyclic module $H^\#_{(\sigma, \delta)}$. Recall that
$(\sigma,\delta)$ is a modular pair in involution if $\sigma$ is a
group--like element in $H$, $\delta:H\to k$ is a morphism of
algebra and the twisted antipode $S_\delta$ is involutive, see
\cite{CM1,CM2}. The cyclic cohomology of $H^\#_{(\sigma, \delta)}$
is called the cyclic cohomology of $H$ and it is denoted
$\mathrm{HC}^\ast_{(\sigma, \delta)}(H)$. One of the features of
$\mathrm{HC}^\ast_{(\sigma, \delta)}(H)$ is that, for a given
algebra $A$ on which $H$ acts and a given $H$--invariant trace
$\tau:A\to k$, there is a canonical morphism:
\[
\gamma^\ast_\tau:\mathrm{HC}^\ast_{(\sigma,
\delta)}(H)\longrightarrow\rmu{HC}{\ast}(A).
\]

Connes--Moscovici cocyclic module was generalized by Khalkhali and
Rangipour in \cite{KR}. Instead of working with modular pairs,
they consider a $\sigma$--compatible Hopf triple $(A,H,M)$, that
is, a Hopf algebra $H$ that coacts to the right on an algebra $A$
together with a left $H$--module and a suitable group--like
element $\sigma\in{H}$. For such a triple $(A,H,M)$ they define a
cyclic module, which in the particular case $A:=H$ and
$M:=k_\delta$ coincides with Connes--Moscovici construction.\vsd

Let us remark that almost all algebras appearing in the above
results are examples of Hopf algebras. Exceptions are only almost
symmetric algebras. Nevertheless, these algebras are
$\Ulie{g}$--Hopf Galois extensions of $k$. Classical Galois
extensions, strongly graded algebras and $H$--crossed products are
other examples of Hopf Galois extension. In spite of the richness
of examples the theory of Hopf Galois extension represents an
unifying setting. By studying Hopf Galois extensions in general,
instead of working with particular examples, the results become
more general and proofs more natural.

Furthermore, dealing with extensions, is more natural to work with
a relative variant of the homology that we are interested in. On
the one hand, in many cases, the relative homology is easier to
compute and, on the other hand, if the subalgebra has ``nice''
homological properties (like separability) then the relative
homology and the usual one are identical.

In this paper we intend to exploit both the unifying character of
Hopf Galois extensions and the simplicity of relative homology.
Our goal is to extend some of the above results and, at the same
time, to show that they are closely related, although they were
obtained using completely different methods.\vsd

Let $A/B$ be an extension of $k$--algebras. We start by recalling
the definition of $\rmd{HH}{\ast}(A/B)$ and $\rmd{HC}{\ast}(A/B)$,
the $B$--relative Hochschild homology of $A$ and, respectively,
the $B$--relative cyclic homology. For shortness we shall call
them the Hochschild, respectively cyclic, homology of the
extension $A/B$. As in the non--relative case (when $B=k$), they
are defined by constructing a certain cyclic object
$\rmd{Z}{\ast}(A/B)$. Then we show that the theory of relative
left derived functors (with respect to a certain projective class
of epimorphisms) can be used to compute $\rmd{HH}{\ast}(A/B)$. As
an immediate consequence it follows that
$\rmd{HC}{\ast}(A/B)\simeq \rmd{HC}{\ast}(A/k)$, whenever $B$ is
separable $k$--algebra.

The properties of Hopf Galois extension that we need are proved in
the second section. By definition an extension $A/B$ is called
Hopf Galois if there is a Hopf algebra $H$ that coacts on $A$ such
that the subalgebra of coinvariant elements is $B$ and a certain
canonical map is bijective. To emphasize the role that $H$ plays
we shall say that $A/B$ is an $H$--Galois extension. Let us denote
the quotient $A/[A,B]$ by $A_B$, where
$[A,B]=\{ab-ba\mid\;a\in{A},\;b\in{B}\}$. It is well-known that
$A_B$ is a left $H$--module and a right $H$--comodule, and one can
check that these structures are compatible:
\begin{eqnarray}
\rho(h\overline{a})&=&\tsum h\de{2}\overline{a}\an{0}\otimes
h\de{3}a\an{1}Sh\de{1},\label{ec:1}\\
\overline{a}&=&\tsum a\an{1}\overline{a}\an{0}.\label{ec:2}
\end{eqnarray}
In these relations $\overline{a}$ denotes the class of $a\in A$ in
$A_B$, and we used the notation $\Delta(h)=\sum h\de{1}\otimes
h\de{2}$ and
$\rho(\overline{a})=\sum\overline{a}\an{0}\otimes\overline{a}\an{1}$.

Using these properties of $A_B$ we construct a new cyclic object
$\rmd{Z}{\ast}(H,A_B)$ that depends only on $H$ and $A_B$. The
main results of this section asserts that $\mathrm{Z}_\ast(A/B)$
and $\rmd{Z}{\ast}(H,A_B)$ are isomorphic cyclic objects. Thus, in
particular, $\rmd{HC}{\ast}(A/B)$ is isomorphic to the cyclic
homology of $\rmd{Z}{\ast}(H,A_B)$.

As we have already noticed, the construction of
$\rmd{Z}{\ast}(H,A_B)$ involves only the fact that $A_B$ is a left
$H$--module and a right $H$--comodule such that (\ref{ec:1}) and
(\ref{ec:2}) hold true. Therefore in the definition of
$\rmd{Z}{\ast}(H,A_B)$ we can replace $A_B$ with an arbitrary left
$H$--module $M$ which is a right $H$--comodule too such that
relations (\ref{ec:1}) and (\ref{ec:2}) hold. We shall call such
an object $M$ a modular crossed module and we shall denote the
corresponding cyclic object by $\rmd{Z}{\ast}(H,M)$. See
Definition \ref{de:CrossedModule} for the definition of modular
crossed modules and Theorem \ref{te:CyclicObject} for the
construction of $\rmd{Z}{\ast}(H,M)$. The Hochschild, respectively
cyclic, homology of $\rmd{Z}{\ast}(H,M)$ will be denoted by
$\rmd{HH}{\ast}(H,M)$, respectively $\rmd{HC}{\ast}(H,M)$.

The properties of the category $\mathcal{CM}_m(H)$ of all modular
crossed module are investigated in Section 4. For example here we
prove that $\mathcal{CM}_m(H)$ is  an abelian category with enough
injective objects. Also, for a Hopf algebra $H$ and a Hopf
subalgebra $K\subseteq H$, we construct a functor
\[
\rmd{Ind}{K}^HM:\mathcal{CM}_m(H)\longrightarrow
H\mathcal{CM}_m(H),\qquad \rmd{Ind}{K}^HM:=H\otimes_HM.
\]
The main goal of this part is the computation of
$\rmd{HC}{\ast}(\rmd{Ind}{K}^HM)$ under the assumptions that $K$
is cocommutative and $M$ is a direct sum of one dimensional
comodules such that the group-like associated to this
decomposition are central in $K$. The result that we obtain is
stated in Theorem \ref{te:General}. The proof is based on a
variant of Shapiro's Lemma, that allows us to reduce the
computation to the calculation of . Since $\rmd{HC}{\ast}(K,-)$
commutes with direct sums we may assume that $\rho(m)=m\otimes g$,
for every $m\in M$, where $g$ is a certain central group--like
element in $K$. The computation of $\rmd{HC}{\ast}(K,M)$,  for
$g=1$, is done in Theorem \ref{th:228}. The case
$\rmd{ord}{}(g)<\infty$ is treated in Corollary \ref{co:232} and
the case $\rmd{ord}{}(g)=\infty$ is considered in Proposition
\ref{pr:030407}.

Let us have a closer look to a modular crossed module $M$ over a
group algebras $H:=kG$. One can easily see that
\[\textstyle
 M\simeq\bigoplus_{g\in{t(G)}}\rmd{Ind}{kG_g}^{kG}M_g,
\]
where $t(G)$ is a coset for the set of conjugacy classes in $G$,
and $G_g$ is the centralizer of $g$ in $G$. Moreover, the comodule
structure of $M_g$ is defined such that $\rho(m)=m\otimes g$, for
every $m\in M_g$. Hence we can apply Theorem \ref{th:228} to
compute the cyclic homology of $kG$ with coefficients in $M$, see
Corollary \ref{co:group-algebras}.

A $G$--strongly graded algebra $A=\bigoplus_{g\in G} A_g$ is a
$kG$--Galois extension of $B:=A_1$. Thus the cyclic homology of
$A/B$ can also be computed using our method. As a matter of fact
this result is obtained by taking $M=A_B$ in Corollary
\ref{co:group-algebras}. By specializing $A$ to $kG$ we give a new
proof of the computation of cyclic homology of group algebras,
performed by Burghelea. In the last part of this section we
compute the cyclic homology of quantum tori, see Theorem
\ref{te:Torus}.

In Section 6, we consider the case of enveloping algebras of Lie
algebras. Let $\Lie{g}$ be a Lie algebra and let $\Ulie{g}$ be its
enveloping algebra. We show that every modular crossed module $M$
can be filtered in a canonical way such that the graded
associated, which obviously is a modular crossed module, has a
trivial comodule structure. By using the computation from Theorem
\ref{th:228} we construct a spectral sequence converging to
$\rmd{HC}{\ast}(\Ulie{g},M)$. Finally, we show that the cyclic
homology of an almost symmetric algebra $A$ is isomorphic to
$\rmd{HC}{\ast}(\Ulie g,M)$, where $\Lie g$ is a certain Lie
algebra and $M$ is a certain modular crossed module over $\Ulie g$
associated to $A$.\vst

\noindent\textbf{Notation}\vsu

Throughout the paper $k$ will denote a field. The tensor product
of two vector spaces will be denoted by $\otimes$. For the
enveloping algebra of an algebra $A$ we whall use the notation
$A^e$. By definition $A^e$ is the algebra $A\otimes A^{opp}$,
where $A^{opp}$ is the opposite algebra structure on $A$.

Let $A$ be an algebra. By definition an $(A,A)$--bimodule is a
left module over $A^e$. The category of $(A,A)$--bimodules will be
denoted by $\ama$.

If $B$ is a subalgebra of $A$ we shall say that $A$ is an algebra
extension of $B$ and we shall write $A/B$.

\section{Hochschild and cyclic homology of extensions of algebras}

\markboth{\sc{Hochschild and cyclic homology of extensions of
algebras}} {\sc{Hochschild and cyclic homology of extensions of
algebras}}

In this section  we shall recall the definition and the basic
properties of relative Hochschild homology and of relative cyclic
homology. We begin by recalling some properties of the cyclic
tensor product, defined by Quillen in \cite{Q}.

\begin{definition}
Suppose that $M$ and $N$ are two $(B,B)$-bimodules. The cyclic tensor
product of $M$ and $N$ is $M{\widehat{\otimes}_B}N:=(M{\otimes_B}%
N)\otimes_{B^e}B$.
\end{definition}

\begin{noname}
For any $(B,B)$--bimodule $X$ we denote by $X_B$ the quotient of
$X$ by the subspace $[X,B]$ generated by all commutators
$[x,b]:=xb-bx$. Since $X\otimes_{B^e}B\simeq{X_B}$, it follows
that
\begin{equation*}
M{\widehat{\otimes}_B}N\simeq(M{{\otimes}_{B}\,}N)/[M{{%
\otimes}_{B}\,}N,B]=(M\otimes_B N)_B.
\end{equation*}
We also have $M{\widehat{\otimes}_B}N\simeq M\otimes_{B^e}N$. The
identification is defined by the canonical map $
m\otimes_{B^e}n\mapsto m{\widehat{\otimes}_B}n$, where
$m{\widehat{\otimes}_B}n:=m{\otimes_B}n+[M{\otimes_B}N,B]$.
\end{noname}
The cyclic tensor product can be defined for an arbitrary, but finite, number of
$(B,B)$--bimodules $M_1$, \ldots, $M_n$ by
\begin{equation*}
M_1{\widehat{\otimes}_B}\cdots{\widehat{\otimes}_B}M_n:=(M_1{\otimes_B}\cdots%
{\otimes_B}M_n)\otimes_{B^e}B.
\end{equation*}
\begin{lemma}\label{le:Maps}
Let $A/B$ be an extension of algebras and let $M_1$, \ldots, $M_n$
be objects in $\bmb$. Then:

a) $(M_1{\otimes_B}\cdots{\otimes_B}M_i){\widehat{\otimes}_B}(M_{i+1}{%
\otimes_B}\cdots{\otimes_B}M_n)\simeq M_1{\widehat{\otimes}_B}\cdots{\widehat{%
\otimes}_B}M_n$.

b) There exists a map $t_{M_1,\ldots,M_n}:M_1{\widehat{\otimes}_B}\cdots{%
\widehat{\otimes}_B}M_n\longrightarrow M_n{\widehat{\otimes}_B}M_1\wotB\cdots{%
\widehat{\otimes}_B}M_{n-1}$ such that $t_{M_1,\ldots,M_n}(m_1{\widehat{\otimes}%
_B}\cdots{\widehat{\otimes}_B}m_n)=m_n{\widehat{\otimes}_B}m_1\wotB\cdots{\widehat{%
\otimes}_B}m_{n-1}$.
\end{lemma}
\begin{proof}
a) Both vector spaces are isomorphic to $(M_1{\otimes}_B\cdots{\otimes}_BM_n)%
{\otimes}_{B^e}B$.

b) Straightforward.
\end{proof}

\begin{noname}\label{nn:simplex}
Suppose now that $A/B$ is an extension of algebras and that $M$ is a $(B,B)$--bimodule.
For every $n\in\mathbb{N}^\ast$ and every $0\leq{i}\leq{n}$ there are
well--defined maps $d_i:M{\widehat{\otimes}_B}A^{{{\otimes}_B}%
n}\longrightarrow{M}{\widehat{\otimes}_B}A^{{{\otimes}_B}(n-1)}$ such that
\begin{equation*}
d_i(m{\widehat{\otimes}_B}{a^1}{\widehat{\otimes}_B}\cdots{\widehat{\otimes}%
_B}a^n)= \left\{
\begin{array}{ll}
ma^1\widehat\otimes{}_B {}a^2\widehat\otimes_B\cdots\widehat\otimes_B%
a^n, & \text{ for }i=0, \\
m{\widehat{\otimes}_B}a^1\widehat\otimes{}_B\cdots{\widehat{\otimes}_B}a^ia^{i+1}{\widehat{\otimes}%
_B}\cdots{\widehat{\otimes}_B}{a^n}, & \text{ for } 0<i<n,\\
a^nm\widehat{\otimes}_B{}a_1\wotB{}\cdots{\widehat{%
\otimes}_B}a^{n-1},& \text{ for } i=n.
\end{array}
\right.
\end{equation*}
If $m:A\otimes_B A\to A$ denotes the multiplication in $A$ and
$\mu_r:M\otimes_B A\to M$ defines the right module structure then
\[
d_0=\mu_r\wotB A^{{\otimes_B}(n-1)} \qquad \text{and} \qquad d_i:=A^{{\otimes_B}(i-1)}
\wotB m\wotB A^{{\wotB} (n-i)}
\]
for every $0<i<n$. Analogously, if $\mu_l$ defines the left module structure on $M$
then we define $d_n:=[\mu_l\wotB A^{{\otimes_B}(n-1)}]\circ{}t_{M,A,\ldots,A}$.

For $0\leq{i}\leq{n}$ one can also define ${}s_i:M{\widehat{\otimes}_B}A^{{{\otimes}_B}%
n}\longrightarrow{M{\widehat{\otimes}_B}A^{{{\otimes}_B}(n+1)}}$ by:
\begin{equation*}
{}s_i(m{\wotB}a^1{\wotB}\cdots{\wotB}a^n) =m\wotB a^1 \wotB \cdots{%
\widehat{\otimes}_B}a^i{\widehat{\otimes}_B}1{\widehat{\otimes}_B}a^{i+1}{%
\widehat{\otimes}_B} \cdots{\widehat{\otimes}_B}a^n.
\end{equation*}
The map $t_n:=t_{A,\ldots,A}$ will play an important role in our paper, since it will
be used to associate a cyclic object to any extension $A/B$ of $k$--algebras.
\end{noname}

\begin{theorem}
\label{thm:simplicial-object} Let $A/B$ be an extension of
$k$--algebras and let $M\in {}_{A}\mathfrak{M}_{A}$.

a) If $\mathrm{Z}_n(A/B,M):=M{\widehat{\otimes}_B}A^{{{\otimes}_B}n}$ then
$\mathrm{Z}_\ast(A/B,M)$ is a simplicial object with face maps $d_0$, \ldots, $d_n$ and
degeneracy maps  $s_0, \ldots,s_{n}$.

b) The simplicial object $\mathrm{Z}_\ast(A/B,A)$ is a cyclic
object with
the cyclic structure defined by $t_n:A^{{\widehat{\otimes}_B}%
(n+1)}\longrightarrow{A^{{\widehat{\otimes}_B}(n+1)}}$, for every $n\geq0$.
\end{theorem}

\begin{proof}
By (\ref{nn:simplex}), $d_\ast$, $s_\ast$ and $t_\ast$ are well defined. One can prove
that they define a simplicial object, respectively a cyclic object, as in the case when
$B=k$, see for example \cite[page 330]{We}.
\end{proof}

\begin{definition}
\label{de:StandardComplex}The \emph{Hochschild homology of an
extension $A/B$ with coefficients in $M$} is the homology of the
complex $(\mathrm{C}_{\ast }(A/B,M),\rmd{b}{\ast})$, associated to
$\mathrm{Z}_\ast(A/B,M)$. The Hochschild homology of $A/B$ will be
denoted by $\rmd{HH}{\ast}(A/B,M)$.

To simplify the notation we shall write $\rmd{Z}{\ast}(A/B)$,
$\rmd{C}{\ast}(A/B)$ and $\rmd{HH}{\ast}(A/B)$ for
$\rmd{Z}{\ast}(A/B,A)$, $\rmd{C}{\ast}(A/B,A)$ and
$\rmd{HH}{\ast}(A/B,A)$, respectively.
\end{definition}

\begin{noname}
Following \cite[Section 9.6]{We} we associate to the cyclic object
$\rmd{Z}{\ast}(A/B)$ a double complex $\ccd(A/B)$ that will be
called the \emph{Tsygan double complex} (we adopt the convention
that in a double complex every square is anti--commutative). Its
columns in even degrees are equal to $\mathrm{C}_{\ast }(A/B),$
while the columns in odd degrees are equal to the acyclic complex
$\mathrm{C}^a_{\ast }(A/B)$. By definition $\mathrm{C}^a_{\ast
}(A/B)=\mathrm{C}_{\ast }(A/B)$ as graded vector spaces. The
differentials of $\mathrm{C}^a_{\ast }(A/B)$ are $-b_\ast'$, where
$b_\ast'=-\sum^{n-1}_{i=0}(-1)^id_i$. The $q^{th}$ row of the
double complex is the periodic complex that computes the homology
of the cyclic group of order $q+1$ acting on $A^{\wotB q+1}$ via
$(-1)^qt_q$.
\end{noname}

\begin{definition}
The \emph{cyclic homology} $\rmd{HC}{\ast}(A/B)$ \emph{of an
extension} $A/B$ is the homology of the total complex of
$\ccd(A/B)$.
\end{definition}
\begin{noname}By \cite[Proposition 9.6.11]{We} it results
immediately that the following sequence is exact.
\[
\cdots\to\rmd{HC}{n+1}(A/B)\overset{S}{\to}
\rmd{HC}{n-1}(A/B)\overset{B}{\to}\rmd{HH}{n}(A/B)\overset{I}{\to}
\rmd{HC}{n}(A/B)\to\cdots
\]
It will be called the $SBI$--\emph{sequence}.
\end{noname}
\begin{noname}\label{cohomology}
We now want to give an equivalent interpretation of Hochschild
homology of an extension $A/B$. More precisely we shall prove that
Hochschild homology of $A/B$ is the left $\calE$--derived functor
of $A\otimes_{A^e}(-):\ama\to\mathcal{M}_k$, where $\calE$ is the
projective class of all epimorphisms of $(A,A)$--bimodules that
splits as morphisms of $(B,B)$--bimodules. To show that $\calE$ is
indeed a projective class of epimorphism and to construct an
$\calE$--projective resolution of $A$ we shall use the formalism
of monoidal categories, that can be found for example in
\cite{MS}.

Algebras, modules, bimodules, etc. can be defined in an arbitrary
monoidal category. All general definitions that we need can be
found for example in \cite{MS}. The category that we work in is
the category of all $(B,B)$--bimodules. It is monoidal with
respect to:
\begin{equation*}
\otimes _{B}:{}_{B}\mathfrak{M}_{B}\times
{}_{B}\mathfrak{M}_{B}\longrightarrow {}_{B}\mathfrak{M}_{B},\quad
(M,N)\mapsto {M\otimes _{B}{N}}.
\end{equation*}
The unit object in $_{B}\mathfrak{M}_{B}$ is $B$, and the
associativity and the unit constraints are the canonical ones. One
can easily check  that to give an algebra in $(_B\mathfrak{M}
_B,\otimes_B,B)$ is equivalent to give a morphism of unitary rings $%
B\to{A}$, see \cite[Example 1.6(d)]{MS}. Thus any extension of
algebras $A/B$ can be thought as an algebra in
$_{B}\mathfrak{M}_{B}$.

Moreover, if $A/B$ is an extension of algebras, that is an algebra
in $_{B}\mathfrak{M}_{B}$, then the category of $(A,A)$--bimodules
in $({}_{B}\mathfrak{M}_{B},\otimes _{B},B)$ is isomorphic to
${}_{A}\mathfrak{M}_{A}$, the category of left modules over the
enveloping algebra $A^{e}$.

By applying \cite[Proposition 1.15]{MS} to our setting it results
that the class of all epimorphism in $\ama$ that splits in $\bmb$
is projective. Moreover, by \cite[Theorem 1.20]{MS} it results
that
\begin{equation*}
\beta _{\ast }(A/B):\qquad 0\longleftarrow {A}\overset{d_{0}}{%
\longleftarrow }A\otimes _{B}{A}\overset{d_{1}}{\longleftarrow }\cdots
\overset{d_{n}}{\longleftarrow }A^{\otimes _{B}(n+2)}\longleftarrow \cdots
\end{equation*}
is an $\mathcal{E}$--relative projective resolution of $A$ in $\ama$. The differentials
of
$\beta _{\ast }(A/B)$ are given by $d_{n}(a^{0}\otimes _{B}\cdots \otimes _{B}{a^{n+1}}%
)=\sum_{i=0}^{n+1}(-1)^{i}a^{0}\otimes _{B}\cdots \otimes _{B}{a^{i}a^{i+1}}%
\otimes _{B}\cdots \otimes _{B}{a^{n+1}}.$

\label{homology} Thus to compute the left $\calE$--derived
functors of $F:=A\otimes_{A^e}(-):\ama\to \calM_k$ we can use the
resolution $\beta _{\ast }(A/B)$. The corresponding complex is
$\mathrm{C}_{\ast }(A/B,M)$ since we can make the following
identifications:
\begin{equation*}
A^{\otimes _{B}(\ast+2)}\otimes _{A^{e}}M\simeq%
A^{\otimes_{B}\ast}\otimes_{B^{e}}A^{e}\otimes_{A^{e}}M\simeq%
A^{\otimes _{B}\ast}\otimes _{B^{e}}M\simeq%
A^{{{\otimes }_{B}\ast}}\wotB M=\mathrm{C}_{\ast}(A/B,M).
\end{equation*}
Furthermore one can see easily that through the identifications above the differential
maps corresponds to those of $\mathrm{C}_{\ast }(A/B,M).$
\end{noname}

\noindent Summarizing we have the following theorem.

\begin{theorem}
\label{thm:Hoch_Rel=Hoch_Mon} Suppose that $A/B$ is an extension
of $k$--algebras and $M$ is an $(A,A)$--bimodule. Then
\[
\rmd{HH}{\ast}(A/B)=\rmd{L}{\ast}^\calE F(A,M),
\]
where $F:=A\otimes_{A^e}(-):\ama\to \calM_k$.
\end{theorem}

\begin{noname}\label{no:separable} Let $A/B$ and $B/C$ be two extensions of
algebras. As an application of this theorem we shall compare the
Hochschild homology groups $\rmd{HH}{\ast}(A/B,M)$ and
$\rmd{HH}{\ast}(A/C,M)$. More precisely we shall prove that these
groups are isomorphic if the extension $B/C$ is separable.

Recall that an extension of algebras $B/C$ is called \emph{separable} if the
multiplication map $B\otimes_C B\rightarrow B$ has a section in $\bmb$. Regarding $B$
as an algebra in the monoidal category $({}_C\mathfrak{M}_C,\otimes_C,C)$ we can check
immediately that $B/C$ is separable iff the algebra $B$ in
$({}_C\mathfrak{M}_C,\otimes_C,C)$ is separable, see \cite{MS}. The main
characterization of separable algebras in a monoidal category, proved in \cite{MS} and
applied to $({}_C\mathfrak{M}_C,\otimes_C,C)$, states that $B/C$ is separable iff any
epimorphism $f:X\rightarrow Y$ in $\bmb$ that splits in ${}_C\mathfrak{M}_C$ has a
section in $\bmb$.
\end{noname}

\begin{proposition}\label{pr:HochSep}
Let $A/B$ and $B/C$ be two extensions of
algebras. If $B/C$ is separable and $M$ is an $(A,A)$--bimodule
then
\[
\rmd{HH}{\ast}(A/C,M)\simeq\rmd{HH}{\ast}(A/B,M).
\]
The isomorphism is induced by the canonical map
$\rmd{C}{\ast}(A/C,M)\rightarrow\rmd{C}{\ast}(A/B,M)$.
\end{proposition}

\begin{proof}
Let $\calE_B$ be the projective class of all epimorphisms in
$\ama$ that have a section in $\bmb$, and let us define $\calE_C$
similarly. Obviously, if $P\in\ama$ is $\calE_C$--projective then
it is $\calE_B$--projective too. Thus any $\calE_C$--projective
resolution of $A$ is made of $\calE_B$--projective bimodules over
$A$. Let $(\Pd,\dd)\longrightarrow A$ be such a resolution. We
claim that $(\Pd,\dd)$ is an $\calE_B$--projective resolution of
$A$. Indeed, let us write $d_\ast$ as a composition
$\mu_\ast\circ\varepsilon_\ast$, where $\mu_\ast$ is an
epimorphism of $(A,A)$--bimodules and $\varepsilon_\ast$ is a
monomorphism in $\ama$. By the definition of $\calE_B$--projective
resolutions, see \cite{HS}, we have to prove that every $\mu_\ast$
has a section of $(B,B)$--bimodules. As $(\Pd,\dd)$ is a
$\calE_C$--projective resolution there is a section of $\mu_\ast$
in $\cmc$. Since $B/C$ is separable it follows by
(\ref{no:separable}) that $\mu_n$ has a section in $\bmb$.

We know that $\rmd{\beta}{\ast}(A/C)$ is an $\calE_C$--projective
resolution of $A$. By the foregoing it is an $\calE_B$--projective
resolution too. On the other hand $\rmd{\beta}{\ast}(A/B)$ is
another $\calE_B$--projective resolution of $A$, and there is a
canonical morphism of resolutions:
\[
\rmd{\beta}{\ast}(A/C)\rightarrow\rmd{\beta}{\ast}(A/B).
\]
that extends the identity of $A$. If we apply to this morphism the
functor $(-)\tp{A^e}M$ we get the canonical map from
$\rmd{C}{\ast}(A/C,M)$ to $\rmd{C}{\ast}(A/B,M)$. We conclude the
proof of the proposition by remarking that this morphism is a
quasi isomorphism of complexes, because two $\calE_B$--projective
resolutions of the same object are homotopically equivalent.
\end{proof}

\begin{corollary} Let $A/B$ be an extension of algebras such that
$B$ is a separable $k$--algebra. Then for every $(A,A)$--bimodule
$M$ we have
\[
\rmd{HH}{\ast}(A/B,M)\simeq\rmd{HH}{\ast}(A,M)\qquad\text{and}\qquad
\rmd{HC}{\ast}(A/B)\simeq\rmd{HC}{\ast}(A),
\]
where $\rmd{HH}{\ast}(A)$ and $\rmd{HC}{\ast}(A)$ are the usual
homology theories (of the extension $A/k$).
\end{corollary}

\begin{proof}The first isomorphism comes directly from the previous proposition, by taking $C=k$.
To prove the second isomorphism one can use the the first part of
the corollary, the $SBI$--sequence and the 5--Lemma.
\end{proof}

\section{$H$--Galois  extensions}
\markboth{\sc{$H$--Galois extensions}} {\sc{$H$--Galois
extensions}}

In this section we recall the definition and basic properties of
an $H$--Galois extension $A/B$. Then we show that the cyclic
homology of $A/B$ can be computed by using a new cyclic object,
that depends only on $H$, the Hopf algebra from the definition of
Galois extensions, and $A/[A,B]$.

Let $A/B$ be an extension of $k$--algebras. For any
$(A,A)$--bimodule $M$ we set
\begin{equation*}
M^{B}:=\{m\in {M}\mid \;bm=mb, \forall b\in {B},\ \forall m\in
{M}\}\quad \text{and}\quad M_{B}:=M/[M,B],
\end{equation*}
where $[M,B]=\langle bm-mb\mid b\in {B},\ m\in {M}\rangle$. Our
aim is to define an algebra structure on $(A{ \otimes _{B}}A)^{B}$
such that $M^{B}$ is a left $(A{\otimes _{B}}A)^{B}$ --module and
$M_{B}$ is a right $(A{\otimes _{B}}A)^{B}$--module. Then we shall
apply this construction to an $H$--Galois extension $A/B$,
rediscovering in that way the Ulbrich--Miyashita action.

\begin{noname}
Suppose that $x,y\in A$. Let
\begin{equation*}
\varphi _{M}^{\prime }:A\otimes _{B}{A}\longrightarrow \mathrm{Hom}%
_{k}(M^{B},M),\;\varphi _{M}^{\prime }(x\otimes _{B}y)(m)=xmy.
\end{equation*}
If $z\in (A{\otimes _{B}}A)^{B}$, then $\varphi _{M}^{\prime }(z)(m)\in {%
M^{B}}$, for any $m\in {M^{B}}$. Hence the restriction of $\varphi
_{M}^{\prime }$ to $(A{\otimes _{B}}A)^{B}$ has the image included into $%
\mathrm{Hom}_{k}(M^{B},M^{B})$. Thus, the restriction of $\varphi
_{M}^{\prime }$ to $(A{\otimes _{B}}A)^{B}$ is a map $\varphi _{M}:(A{%
\otimes _{B}}A)^{B}\longrightarrow \mathrm{End}_{k}(M^{B})$ such that for any
$z=\sum_{i=1}^{n}x_{i}{\otimes _{B}}y_{i}$ in $(A{\otimes _{B}}A)^{B}$ we have:
\begin{equation*}
\varphi _{M}(z)(m)=\tsum\, {x_{i}my_{i}},\quad \forall m\in
{M^{B}}.
\end{equation*}
If we take $M:=A{\otimes _{B}}A$, with its natural structure of $(A,A)$%
--bimodule, then $\varphi _{A{\otimes _{B}}A}$ defines a multiplication
\begin{equation*}
``\cdot ":(A{\otimes _{B}}A)^{B}\otimes (A{\otimes _{B}}A)^{B}%
\longrightarrow (A{\otimes _{B}}A)^{B},\;z\cdot {z^{\prime }}=\varphi _{z}(z^{\prime
}).
\end{equation*}
With respect to this operation $(A{\otimes _{B}}A)^{B}$ becomes an
associative $k$--algebra, with unit $1{\otimes _{B}}1$. It
generalizes the usual enveloping algebra of a $k$--algebra $A$, so
it will be called the \textit{enveloping algebra of} $A/B$.

Note that $\varphi_M$ defines a left $(A{\otimes_B}A)^B$--module
structure on $M^B$ by:
\begin{equation*}
z\cdot{m}=\varphi_M(z)(m),\quad \forall z\in(A{\otimes_B}A)^B\quad
\forall m\in{M^B}.
\end{equation*}
\end{noname}

\begin{noname}
Keeping the notation from the previous paragraph, one can check
easily that the map
\begin{equation*}
\psi^{\prime}_M:A{\otimes_B}A\longrightarrow\mathrm{Hom}_k(M,M_B), \;
\psi^{\prime}_M(x{\otimes_B}y)(m)=ymx+[M,B]
\end{equation*}
is well defined. Also $\psi^{\prime}_M(x)(mb-bm)=0$ for any $z\in(A{\otimes_B%
}A)^B$, $m\in{M}$ and $b\in{B}$, as $\psi^{\prime}_M(z)$ does not
depend on a particular choice of elements $x_1$, \ldots, $x_n$,
$y_1$, \ldots, $y_n$ such that $z=\sum_{i=1}^nx_i{\otimes_B}y_i$.
Hence we obtain a map
\begin{equation*}
\psi_M:(A{\otimes_B}A)^B\longrightarrow\mathrm{Hom}_k(M_B,M_B), \;
\psi_M(\tsum_{i=1}^nx_i{\otimes_B}y_i)(\overline{m})=\tsum_{i=1}^n\overline{%
y_imx_i}
\end{equation*}
that defines a right $(A{\underset{B}{\otimes}\,}A)^B$--module structure on $%
M_B$ by
\begin{equation*}
\overline{m}\cdot{z}=\psi_M(z)(\overline{m}).
\end{equation*}
\end{noname}

\begin{noname}\label{pa:Galois}
Suppose now that $A/B$ is an $H$--Galois extension, where $H$ is a
given Hopf algebra over $k$. This means that $A$ is a right
$H$--comodule via an algebra map
$\rho:A\longrightarrow{A\otimes{H}}$, such that $B=\{a\in{A}
\mid\;\rho(a)=a\otimes1\}$ and the canonical map
\begin{equation*}
\beta:A{\otimes_B}A\longrightarrow{A\otimes{H}}, \;
\beta(a{\otimes_B} y)=a\rho(y)
\end{equation*}
is bijective. Using Sweedler's notation, $\rho(y)=\tsum
y\an{0}\otimes y\an{1}$, we have
\begin{equation*}
\beta(x{\otimes_B}y)=\tsum{xy_{\langle{0}\rangle}}\otimes{y_{\langle{1}
\rangle}}.
\end{equation*}
Let us remark that $\beta$ is a morphism of $(B,B)$--bimodules
where $A\otimes{H}$ is regarded as a bimodule with the structure
induced from that one of $A$. Hence $\beta$ induces a $k$--linear
isomorphism $
\overline{\beta}:(A{\otimes_B}A)^B\longrightarrow{A^B\otimes{H}}$,
we can define a $k$--linear map
$\kappa:H\longrightarrow(A{\otimes_B}A)^B$ by $
\kappa:=(\overline{\beta})^{-1}\circ{i}$. Here
$i:H\longrightarrow{A^B\otimes{H}} $ denotes the canonical map
$i(h)=1\otimes{h}$. For every $h\in{H}$ we will use the notation
$\kappa(h)=\sum\kappa^1(h){\otimes_B}\kappa^2(h)$, therefore, by
the definition of $\overline{\beta}$ and $\kappa$, we have:
\begin{equation*}
\tsum\kappa^1(h)\kappa^2(h)_{\langle{0}\rangle}\otimes\kappa^2(h)_{\langle{1}%
\rangle}=1\otimes{h}.
\end{equation*}
Thus, for $h$, $k\in{H}$ we get:
\begin{equation*}
\overline{\beta}(\kappa(h)\kappa(k))
=\tsum\kappa^1(h)\kappa^1(k)\kappa^2(k)_{\langle{0}\rangle}\kappa^2(h)_{\langle{0}\rangle}\otimes%
\kappa^2(k)_{\langle{1}\rangle}\kappa^2(h)_{\langle{1}\rangle}
=1\otimes{kh} =\overline{\beta}(\kappa(kh)).
\end{equation*}
It result that $\kappa$ is an anti--morphism of algebras. Summarizing, we get the
following proposition.
\end{noname}

\begin{proposition}
Let $A/B$ be an $H$--Galois extension. Suppose that $M$ is an $(A,A)$%
--bimodule. Then:

\indent a) $M^B$ is a right $H$--module with the structure
\begin{equation*}
m\cdot{h}=\tsum\kappa^1(h)m\kappa^2(h), \forall h\in{H},\quad
\forall m\in{M^B}.
\end{equation*}
\indent b) $M_B$ is a left $H$--module with the structure
\begin{equation*}
h\cdot\overline{m}=\tsum\overline{\kappa^2(h)m\kappa^1(h)},
\forall h\in{H},\quad \forall m\in{M_B}.
\end{equation*}
Both structures are functorial in $M$.
\end{proposition}

\begin{proof}
$M^B$ is a left $(A{\otimes_B}A)^B$--module. Since $\kappa$ is an anti--morphism of
algebras it defines a right $H$--module structure on $M^B$. Similarly we can define the
left action of $H$ on $M_B$.
\end{proof}

\begin{remark}
Both structures have already appeared in \cite{St}, where they are called the
Ulbrich--Miyashita actions.
\end{remark}

\noindent Some other useful properties of $\kappa$ are listed and
proved in the next proposition.

\begin{proposition}\label{pr:Galois}
Let $A/B$ be an $H$--Galois extension. If $M$  is an $(A,A)$--bimodule, $h\in H$,
$a,x\in A$ and $m\in M$ then the following relations hold true.
\begin{equation}\label{ec:tau1}
\tsum\kappa^1(h_{({1})}){\wotB\,}\kappa^2(h_{({1})})\otimes{h_{({2%
})}} =\tsum\kappa^1(h){\wotB\,}\kappa^2(h)_{\langle{0}%
\rangle}\otimes\kappa^2(h)_{\langle{1}\rangle}.
\end{equation}
\begin{equation}\label{ec:tau2}
\tsum\kappa^1(h_{({2})}){\wotB\,}\kappa^2(h_{({2})})\otimes{Sh_{({%
1})}} =\tsum\kappa^1(h)_{\langle{0}\rangle}{\wotB\,}%
\kappa^2(h)\otimes\kappa^1(h)_{\langle{1}\rangle}.
\end{equation}
\begin{equation}\\
\tsum\kappa^1(h)\kappa^2(h)=\varepsilon(h)\label{ec:tau3}.
\end{equation}
\begin{equation}
\rho(h\overline{a})=\tsum h\de{2}\overline{a}\an{0}\otimes
h\de{3}a\an{1}Sh\de{1}.\label{ec:tau5}
\end{equation}
\begin{equation}
\tsum a\an{1}\overline{ma}\an{0}=\overline{am}.\label{ec:tau4}
\end{equation}
\end{proposition}

\begin{proof} All assertions are consequences of
the equations (a)--(g) in \cite[Remark 3.4.2]{Sch}. Indeed, in
\cite{Sch} the notation used for $\beta^{-1}(1\otimes h)$ is $\sum
r_i (h)\otB l_i(h)$, so
\[
\tsum r_i (h)\otB
l_i(h)=\kappa(h)=\tsum\kappa^{1}(h)\otB\kappa^{2}(h).
\]
Thus we can substitute the element $\sum r_i (h)\otB l_i(h)$ by
$\sum\kappa^{1}(h)\otB\kappa^{2}(h)$ in \cite[Relations (d) and (e)]{Sch}. We get
(\ref{ec:tau1}) and (\ref{ec:tau2}) by applying $\pi\otimes H$ to the equations that we
obtain, where $\pi:A\otB A\longrightarrow A\wotB A$ denotes the canonical projection.
Relation (c) in \cite[Remark 3.4.2]{Sch} is exactly (\ref{ec:tau3}) if we use our
notation for $\kappa(h)$.

Relations (\ref{ec:tau1}) and (\ref{ec:tau2}) together imply:
\[
\tsum\kappa^1(h_{({2})}){\wotB\,}\kappa^2(h_{({2})})\otimes{Sh_{({%
1})}}\otimes h\de{3} =\tsum\kappa^1(h)_{\langle{0}\rangle}{\wotB\,}%
\kappa^2(h)\an{0}\otimes\kappa^1(h)_{\langle{1}\rangle}
\otimes\kappa^2(h)_{\langle{1}\rangle}.
\]
Therefore we obtain:
\[
\tsum\overline{\kappa^2(h_{({2})})a\an{0}\kappa^1(h_{({2})})}\otimes h\de{3}a\an{1}{Sh_{({%
1})}}=\tsum\overline{\kappa^2(h)_{\langle{0}\rangle}a\an{0}%
\kappa^1(h)\an{0}}\otimes\kappa^2(h)_{\langle{1}\rangle}a\an{1}
\kappa^1(h)_{\langle{1}\rangle},
\]
so (\ref{ec:tau5}) holds true by the definition of the $H$--module
structure on $A_B$ and the fact that $A$ is an $H$--comodule
algebra. By \cite[Remark 3.4.2(b)]{Sch} and the definition of the
cyclic tensor product we have
\[
\tsum ma\an{0}\kappa^{1}(a\an{1})\wotB\kappa^{2}(a\an{1})=m\wotB
a.
\]
By applying $d_2$ (see (\ref{nn:simplex}) for the definition od
$d_2$)  to both sides of this equation and using the definition of
the left action of $H$ on $M_B$, we get (\ref{ec:tau4}).
\end{proof}

\begin{noname}\label{nn:Galois}
The canonical isomorphism
$\beta:A{\otimes_B}A\longrightarrow{A\otimes{H}}$ can be
inductively extended to obtain isomorphisms
$\beta_n:A^{{\otimes_B}(n+1)}\longrightarrow{A\otimes{H^{\otimes{n}}}}$.

To simplify the notation, a tensor monomial
$h^1\otimes\cdots\otimes h^n\in H^{\otimes n}$ will be denoted by
$(h^1,\ldots,h^n)$.

Obviously $\beta _{n}$ is left $A$--linear with respect to the canonical $A$%
--module structures on $A^{{\otimes _{B}}(n+1)}$ and $A\otimes
{H^{\otimes (n+1)}}$. Therefore, for each $(A,A)$--bimodule $M$,
the map $\beta _{n}(M):=M\otimes _{A}\beta _{n}$ is an isomorphism
of $(A,B)$--bimodules.

For $x=m\otB {a^{1}}\otB \cdots \otB {a^{n}}$ we have
\begin{equation}
 \beta _{n}(M)(x)=
 \tsum {m}a_{\langle {0}\rangle }^{1}\cdots {a_{\langle {0}\rangle}^{n}}
 \otimes (a_{\langle {1}\rangle }^{1}\cdots {a_{\langle{1}\rangle }^{n}},
 \ldots,
 a_{\langle {n-1}\rangle }^{n-1}a_{\langle{n}\rangle }^{n},
 a_{\langle {n}\rangle }^{n})
\label{eq:betaM}
\end{equation}
and, if $y=m\otimes({h^{1}},\cdots,{h^{n}}$, then
\begin{equation}\label{eq:betaMInv}
\beta _{n}(M)^{-1}(y)=\tsum {m}\kappa
^{1}(h^{1}){\underset{B}{\otimes }\,}\kappa ^{2}(h^{1})\kappa
^{1}(h^{2}){\underset{B}{\otimes }\,}\cdots {\underset{B}{\otimes
}\,}\kappa ^{2}(h^{n-1})\kappa ^{1}(h^{n}){\underset{B}{\otimes
}\,}\kappa ^{2}(h^{n}).
\end{equation}
One can see easily that $\beta _{n}(M)$ is a map of
$(B,B)$--bimodules,
therefore it factorizes to a map: $\overline{\beta }_{n}(M):(M{\otimes _{B}}%
A^{{{\otB}}\/n})_{B}\longrightarrow {M_{B}}\otimes {H^{\otimes
{n}}}$ such that
\begin{equation*}
\begin{xy} \xymatrix{ M\otB{A}^{{\otB}\/n} \ar[d]_{\pi} \ar[rr]^{\beta_n(M)}
&&M\otimes{H}^{\otimes{n}} \ar[d]^{\pi_M\otimes{H}^{\otimes{n}}}\\
(M\otB{A}^{\otB{n}})_B \ar[rr]_{\ \overline{\beta}_n(M)}
&&M_B\otimes{H}^{\otimes{n}} }\end{xy}
\end{equation*}
is commutative. In this diagram $\pi $ and $\pi _{M}$ are the
canonical projections.
\end{noname}

\begin{example}
Let $H$ be a Hopf algebra. If $H$ is regarded as an $H$--comodule via $%
\Delta $, the comultiplication of $H$, then $H/k$ is an
$H$--Galois
extension. The maps $\beta_n:H^{\otimes{n+1}}\longrightarrow{H^{\otimes{n+1}}%
}$ are given in this case by
\begin{equation*}
\beta_n(h^0,\ldots,{h^n})=\tsum({h^0}h^1_{({1})}h^2_{({1})}\cdots{%
h^n_{({1})}},{h^1_{({2})}}\cdots{h^n_{({2})}},\ldots,
h^{n-1}_{({n})}h^n_{({n})},{h^n_{({n+1})}}),
\end{equation*}
and their inverses satisfy:
\begin{equation*}
\beta^{-1}_n(h^0,{h^1},\ldots\,{h^n})=\tsum({h^0}Sh^1_{({1}%
)},{h^1_{({2})}}Sh^2_{({1})},\ldots,{h^{n-1}_{({2})}}%
Sh^n_{({1})},{h^n_{({2})}}).
\end{equation*}
The formula for $\beta^{-1}_n$ has been deduced by using
\begin{equation*}
\beta^{-1}(1\otimes{h})=\tsum Sh_{({1})}\otimes{h_{({2})}}.
\end{equation*}
\end{example}

\begin{noname}
$H^{\otimes(n+1)}$ is a right $H$--module with the following
structure:
\[
 (h^1,\ldots,{h^{n+1}})h=
 \tsum (h^1h\de{1},\ldots,{h^{n+1}h\de{n+1}}).
\]
In that case we say that $H^{\otimes(n+1)}$ is a right module via
the diagonal action. Actually this module is a right $H$--comodule
too, with respect to $H^{\otimes n}\otimes \Delta$. In fact,
$H^{\otimes(n+1)}$ is a Hopf module.

We recall that a right $H$--module $M$, together a map $\rho:M\to
M\otimes H$ that defines a comodule structure on $M$, is called a
right Hopf module if
\[
\rho(mh)=\tsum{m}\an{0}h\de{1}\otimes{m}\an{1}h\de{2},\quad\forall
m\in{M},\ \forall h\in{H}.
\]
The structure theorem for Hopf modules says that
$m\otimes{h}\mapsto mh$ defines an isomorphism of Hopf modules
$M^{co(H)}\otimes{H}\simeq M$, where $M^{co(H)}=\{m\in
M\mid\rho(m)=m\otimes1\}$ and $M^{co(H)}\otimes{H}$ is regarded as
a Hopf module with the canonical structure given by the
multiplication and comultiplication of $H$. The inverse map is
\[
 \varphi_M:M\longrightarrow{M^{co(H)}}\otimes{H},\quad
 \varphi_M(m)=\tsum{m}\an{0}\,Sm\an{1}\otimes{m\an{2}}.
\]
For details, see \cite[Theorem~4.1.1]{Sw}.
\end{noname}

\begin{lemma}
If $H$ is a Hopf algebra over a field $k$. For every natural
number $n>0$ there exists an isomorphism of right $H$--modules
$\varphi_n:H^{\otimes{n}}\otimes{H}\longrightarrow
H^{\otimes(n+1)}$, where  $H^{\otimes{n}}\otimes{H}$ has the
canonical $H$--module structure and $H^{\otimes(n+1)}$ is an
$H$--module via $\Delta_n:H\longrightarrow H^{\otimes(n+1)}$, the
$n^{th}$--iterated comultiplication on $H$.
\end{lemma}
\begin{proof}
As we noted before, $M=H^{\otimes(n+1)}$ is a Hopf module.
Moreover, we have $M^{co(H)}=H^{\otimes{n}}$, as the comodule
structure of $M$ is obtained by applying $\Delta$ on the last
factor of $H^{\otimes(n+1)}$. Hence, by the structure theorem of
Hopf modules, we can take $\varphi_n:=\varphi_M^{-1}$.
\end{proof}
\begin{remark}\label{le:030115}
For future references we give the explicit formula for
$\varphi_{n}$ and its inverse. One can see easily  that:
\[
\begin{array}{ll}\varphi_{n}((h^1,\ldots,{h^{n}})\otimes{h})=
 \sum({h^1}h\de{1},\ldots,{h^{n}h\de{n}},{h}\de{n+1})\\
 \varphi_{n}^{-1}(h^1,\ldots,{h^{n+1}})=
 \sum(h^1\,Sh^{n+1}\de{n},\ldots,{h^{n}
 \,Sh^{n+1}\de{1})}\otimes{h}^{n+1}\de{n+1}.
\end{array}
\]
\end{remark}

\section{Cyclic homology of $H$--Galois  extensions}
\markboth{\sc{Cyclic homology of $H$--Galois extensions}}
{\sc{Cyclic homology of $H$--Galois extensions}}

\begin{proposition}\label{pr:C_*(H)}
Let $\wrmd{Z}{\ast}(H):=H^{\otimes(\ast+1)}$.

a) $\wrmd{Z}{\ast}(H)$ is a cyclic object with respect to the following operators (in
degree $n$)
\[
\begin{array}{lll}
 \partial_i=H^{\otimes i}\otimes\varepsilon\otimes H^{\otimes n-i}
 &\qquad\mbox{\em (face maps)},\\
 \mu_i=H^{\otimes i}\otimes\Delta\otimes H^{\otimes n-i}
 &\qquad\mbox{\em (degeneracy maps)},\\
 t_n(h^0,\ldots,h^{n})=(h^{n},h^0,\ldots,h^{n-1})
 &\qquad\mbox{\em (cyclic operator)}.
\end{array}
\]

b) If $\wrmd{C}{\ast}(H)$ is the complex associated to
$\wrmd{Z}{\ast}(H)$ then
$\wrmd{C}{\ast}(H)\overset{\varepsilon}{\longrightarrow}k$ is a
free resolution of $k$ as a right $H$--module.

c) The Hochschild homology $\wrmd{HH}{\ast}(H)$ and cyclic homology
$\wrmd{HC}{\ast}(H)$ of $\wrmd{Z}{\ast}(H)$ are given by
\begin{equation*}
\begin{array}{llll}
\wrmd{HH}{0}(H)=k & \text{and} & \wrmd{HH}{n}(H)=0,& \text{for }
n>0.\\
\wrmd{HC}{2n}(H)=k &\text{and} & \wrmd{HC}{2n+1}(H)=0,&\text{for
}n \geq 0.
\end{array}
\end{equation*}

d) $\wrmd{Z}{\ast}(H)$ is a simplicial object in the category of right $H$--modules. If
$H$ is cocommutative then $t_\ast$ is also a morphism of $H$--modules, so
$\wrmd{Z}{\ast}(H)$ is a cyclic object in $\mm_H$.
\end{proposition}
\begin{proof}
a) Obviously $\partial_i\partial_j=\partial_{j-1}\partial_i$ and
$\mu_i\mu_j=\mu_{j+1}\mu_i$ if $i<j$. The equality $\mu_i\mu_i=\mu_{i+1}\mu_i$ comes
from the fact that $\Delta$ is coassociative. Also
$\partial_i\mu_j=\mu_{j-1}\partial_i$, for $i<j$, and
$\partial_i\mu_j=\mu_{j}\partial_{i-1}$, for $i>j+1$, are straightforward. Finally, for
$i=j$ or $i=j+1$, we have $\partial_i\mu_j=\mathrm{Id}$ as a consequence of the fact
that $\varepsilon$ is the counit of $H$. To prove that $t_n$ is a cyclic operator we
have to check, for $0<i\leq n$, the following relations
\[
\partial_i t_n=t_{n-1}\partial_{i-1}\qquad \mu_i
t_n=t_{n+1}\mu_{i-1}\qquad \partial_0 t_n=\partial_n\qquad \mu_0
t_n=t_{n+1}^2\mu_n.
\]
For example, by evaluating both sides of the first equality at
$(h^0,\ldots,h^{n})$ we obtain
$\varepsilon(h^{i-1})(h^{n},h^0,\ldots,\widehat{h^{i-1}},\ldots,h^{n-1})$.
Similarly one can prove the other relations.

b) We know that every $H^{\tp{} (n+1)}$ is a Hopf module, so by
\cite[Theorem~4.1.1]{Sw} is free.

For $n>0$ let $\widetilde{d}_n:\wrmd{C}{n}(H)\longrightarrow\wrmd{C}{n-1}(H)$ be the
differential of $\wrmd{C}{\ast}(H)$. Take $\widetilde{d}_0=\varepsilon$. Then for every
$n\in\mathbb{N}$ and $x\in\wrmd{C}{n}(H)$ we have
\[
 \widetilde{d}_{n+1}(1\otimes{x})=x-\widetilde{d}_n(x).
\]
So this complex is acyclic and, of course, $\widetilde{d}_0$ is
surjective.

c) By (b) we get $\wrmd{HH}{n}(H)=0$ and that the image of
$\widetilde{d}_1:\rmd{C}{1}(H)\longrightarrow\rmd{C}{0}(H)$ is
$H_+:=\mathrm{Ker}\;\varepsilon$. Hence
$\wrmd{HH}{0}(H)=H/\mathrm{Ker}\,\varepsilon=k$. Let $\wccd(H)$ be the Tsygan double
complex of $\wrmd{Z}{\ast}(H)$. Since all its columns are acyclic, the second spectral
sequence of $\wccd(H)$ degenerates to give us the formulas for $\wrmd{HC}{\ast}(H)$.

d) Since we always consider the right diagonal action on $H^{\otimes(n+1)}$ it is easy
to see that both $\partial_i$ and $\mu_i$ are right $H$--linear, for every $0\leq i\leq
n$. It remains to show that $t_n$ is right $H$--linear if $H$ is cocommutative. But
\[
t_n\left((h^0,\cdots, h^{n})h\right )=\tsum(h^{n}h\de{n+1},h^0h\de{1},\cdots,
h^{n-1}h\de{n}).
\]
As $\tsum(h\de{n+1},h\de{1},\cdots, h\de{n})=\tsum(h\de{1},\cdots,
h\de{n+1})$ we conclude.
\end{proof}

\begin{corollary}\label{co:S(M)}
Let $M$ be a left $H$--module. Then $\rmd{Z}{\ast}(H,M):=\wrmd{Z}{\ast}(H)\otimes_H M$
is a simplicial object. In degree $n$ its face and degeneracy maps are:
\begin{eqnarray}
 \delta_i((h^0,\ldots,h^{n})\otimes_Hm)
 &=&\varepsilon(h^{i})(h^0,\ldots,\widehat{h^{i}},\ldots,h^{n})
 \otimes_Hm,\label{eq:delta_i}\\
 \sigma_i((h^0,\ldots,h^{n})\otimes_Hm)
 &=&(h^0,\ldots,\Delta{h^{i}},\ldots,h^{n})\otimes_Hm.\label{eq:sigma_i}
\end{eqnarray}
The homology of $\rmd{C}{\ast}(H,M)$, the complex associated to $\rmd{Z}{\ast}(H,M)$,
is $\Tor{H}{\ast}(k,M)$.
\end{corollary}
\begin{proof}
The covariant functor $(-)\tp{H} M:\mm_H\longrightarrow \mm_k$ maps a simplicial object
to a simplicial object. Since
$\wrmd{C}{\ast}(H)\overset{\varepsilon}{\longrightarrow}k$ is a free resolution of $k$
we can compute $\Tor{H}{\ast}(k,M)$ as the homology of
$\rmd{C}{\ast}(H,M)=\wrmd{C}{\ast}(H)\otimes_H M$.
\end{proof}
\begin{proposition}\label{pr:IzoSymp}
For any $H$--Galois extension $A/B$ and any $(A,A)$--bimodule $M$ the simplicial
objects $\rmd{Z}{\ast}(A/B,M)$ and $\rmd{Z}{\ast}(H,M_B)$ are isomorphic and
\[
\rmd{HH}{\ast}(A/B,M)\simeq \rmd{H}{\ast}(\rmd{C}{\ast}(H,M_B))\simeq
\rmd{Tor}{\ast}^H(k,M_B).\]
\end{proposition}

\begin{proof}
Let us consider the following sequence of isomorphisms
\begin{equation*}\label{eq:030115}
\rmd{Z}{n}(A/B,M)\stackrel{\overline{\beta}_n(M)}{\longrightarrow}M_B\otimes{H^{\otimes{n}}}
 \stackrel{\simeq}{\longrightarrow}
 H^{\otimes{n}}\otimes{M_B}
 \stackrel{\simeq}{\longrightarrow}
 (H^{\otimes{n}}\otimes{H})\otimes_HM_B
 \stackrel{\simeq}{\longrightarrow}\rmd{Z}{n}(H,M_B).
\end{equation*}
The second morphism is the canonical flip map, the third morphism
is defined by:
$$
    (h^1,\ldots,h^n)\otimes_H\overline{m}
    \mapsto(h^1,\ldots,h^n,1)\otimes_H\overline{m},
$$
and the last one is $\varphi_n\otimes_HM_B$, where $\varphi_{n}$
is explicitly described in Remark~\ref{le:030115}. Let
$\lambda_{n}:\rmd{Z}{n}(A/B,M)\longrightarrow \rmd{Z}{n}(H,M_B)$
be the composition of the above isomorphisms. Then:
\[
 \lambda_n(m\dwotB a^1\dwotB\cdots \dwotB a^n)=\tsum
(\pr{1}{n}{1},\pr{2}{n}{2},\ldots,
a^{n-1}\an{n}a^n\an{n-1}, a^n\an{n},1)\otimes_H%
\overline{{m}\pr{1}{n}{0}}.
\]
If we prove that:
\begin{eqnarray}
 \delta_i\lambda_n(m\wotB a^1\wotB\cdots\wotB a^n)=
 \lambda_{n-1}d_i(m\wotB a^1\wotB\cdots\wotB a^n)\label{eq:d_i}\\
 \sigma_i\lambda_n(m\wotB a^1\wotB\cdots\wotB a^n)
 =\lambda_{n+1}{}s_i(m\wotB a^1\wotB\cdots\wotB a^n)\label{eq:s_i}
\end{eqnarray}
then $\lambda_\ast$ is an isomorphism of simplicial objects.

We first prove  (\ref{eq:d_i}) for $i=0$. If
$A_l=\delta_0\lambda_n(m\wotB a^1\wotB\cdots\wotB a^n)$ we have:
\begin{eqnarray*}
 A_l&=&\tsum\; \varepsilon(\pr{1}{n}{1})(\pr{2}{n}{2},\cdots, a^n\an{n},1)
    \otimes_H \overline{{m}\pr{1}{n}{0}}\\
 &=&\tsum(\pr{2}{n}{1},\cdots,a^n\an{n-1},1)
 \otimes_H\overline{{ma^1}\pr{2}{n}{0}},
\end{eqnarray*}
where for the last equality we used the fact that $\varepsilon$ is
the counit of $H$.

 On the other hand $d_0(m\wotB{a^1}\wotB\cdots\wotB
a^n)=ma^1\wotB{a^2}\wotB\cdots\wotB a^n$, so if we denote by $A_r$ the right hand side
of (\ref{eq:d_i}) we get:
\begin{eqnarray*}
 A_r&=&\lambda_n(ma^1\wotB{a^2}\wotB\cdots\wotB a^n)\\
 &=&\tsum (\pr{2}{n}{1},\pr{3}{n}{2},\ldots,
 a^n\an{n-1},1)
    \otimes_H \overline{{ma^1}\pr{2}{n}{0}}.
\end{eqnarray*}
Thus $A_l=A_r$, so (\ref{eq:d_i}) holds for $i=0$. Similarly one
can prove (\ref{eq:d_i}) for $0<i<n$. For $i=n$ let as denote by
$B_l$ and $B_r$ the left, respectively right, hand sides of
(\ref{eq:d_i}). It is easy to see that
\[
B_r=\tsum (\pr{1}{n-1}{1},\pr{2}{n-1}{2},\ldots,
a^{n-1}\an{n-1},1)\otimes_H \overline{{a^n m}\pr{1}{n-1}{0}}.
\]
On the other hand:
\begin{eqnarray*}
B_l &=&\tsum (\pr{1}{n}{1},\pr{2}{n}{2},\ldots,
    \pr{n-1}{n}{n-1}, a^n\an{n})\otimes_H%
    \overline{{m}\pr{1}{n}{0}})\\
    &=&\tsum (\pr{1}{n-1}{1},\pr{2}{n-1}{2},\ldots,
    \pr{n-1}{n}{n-1}, 1)a^n\an{1}\otimes_H%
    \overline{{m}\pr{1}{n}{0}})\\
    &=&\tsum (\pr{1}{n-1}{1},\pr{2}{n-1}{2},\ldots,
    \pr{n-1}{n}{n-1}, 1)\otimes_H%
    a^n\an{1}\overline{{m}\pr{1}{n}{0}})
\end{eqnarray*}
By (\ref{ec:tau4}) we deduce that $B_l=B_r$.

To finish the proof it remains to show that (\ref{eq:s_i}) holds
true. We shall prove this relation only for $i=0$, as for
arbitrary $0<i\leq n$ one can proceed analogously. Let $C_l$ and
$C_r$ be the left and right hand sides of (\ref{eq:s_i}). It
follows:
\begin{eqnarray*}\textstyle
C_l &=&\sigma_0(\tsum (\pr{1}{n}{1},\pr{2}{n}{2},\ldots,
    \pr{n-1}{n}{n-1}, a^n\an{n},1)\otimes_H%
    \overline{{m}\pr{1}{n}{0}})\\
&=&\tsum (\Delta\left(\pr{1}{n}{1}\right),\pr{2}{n}{2},\ldots,
    \pr{n-1}{n}{n-1}, a^n\an{n},1)\otimes_H%
    \overline{{m}\pr{1}{n}{0}}\\
&=&\tsum (\pr{1}{n}{1},\pr{1}{n}{2},\pr{2}{n}{3},\ldots,
    \pr{n-1}{n}{n}, a^n\an{n+1},1)\otimes_H%
    \overline{{m}\pr{1}{n}{0}}\\
&=& \lambda_{n+1}(m\wotB 1\wotB a^1\wotB\cdots\wotB a^n)%
\end{eqnarray*}
Hence $C_l=C_r$, equality that completes the proof of the fact that $\lambda_\ast$ is
an isomorphism of simplicial objects. By the first part of the proposition it follows
$\rmd{HH}{\ast}(A/B,M)$ and $\rmd{H}{\ast}(\rmd{C}{\ast}(H,M_B))$ are isomorphic graded
vector spaces. The other isomorphism follows by Corollary \ref{co:S(M)}.
\end{proof}

\begin{remark}\label{re:Cat_C}
Let $\mathfrak{C}$ be the category of all quadruples $(A,B,H,M)$, where $A/B$ is an
$H$--Galois extension and $M\in\ama$. The morphisms in $\mathfrak{C}$ are triples
$(f,g,u)$ with $f:A\rightarrow A'$ a morphism of algebras, $g:H\rightarrow H'$ a
morphism of Hopf algebras and $u:M\rightarrow M'$ a morphism of $(A,A)$--bimodules such
that
\begin{eqnarray*}
&\rho(f(a))=\tsum f(a\an{0})\tp{} g(a\an{1}),\\
&u(am)=f(a)u(m)\qquad\text{and}\qquad u(ma)=u(m)f(a).
\end{eqnarray*}
One can see easily that $M\mapsto\rmd{C}{\ast}(A/B,M)$ and
$M\mapsto\rmd{C}{\ast}(H,M_B)$ define two functors from
$\mathfrak{C}$ to the category of chain complexes. Moreover, for
an arbitrary morphism $(f,g,u):(A,B,H,M)\to (A',B',H',M')$ in
$\mathfrak{C}$, the following diagram is commutative
\[
\begin{xy}
\xymatrix{
  \rmd{C}{\ast}(A/B,M)\ar[d]_{\lambda_\ast}\ar[r]
                & \rmd{C}{\ast}(A'/B',M') \ar[d]^{\lambda_\ast}  \\
  \rmd{C}{\ast}(H,M_B) \ar[r]
                & \rmd{C}{\ast}(H',M'_{B'})}
\end{xy}
\]
so
$\lambda_\ast:\rmd{C}{\ast}(A/B,M)\rightarrow\rmd{C}{\ast}(H,M_B)$
is a natural map.
\end{remark}

\begin{corollary}\label{co:LambdaFunct}
Assume that $(A,C,H,M)$ and $(A,B,K,M)$ are two objects in $\mathfrak{C}$. Let
$g:H\rightarrow K$ be a Hopf algebra map such that $(\Id_A,g,\Id_M)$ is a morphism in
$\mathfrak{C}$. If $B/C$ is separable, then the canonical map
$\rmd{C}{\ast}(H,M_C)\rightarrow\rmd{C}{\ast}(K,M_B)$ is a quasi isomorphism of
complexes.
\end{corollary}
\begin{proof}
Apply Proposition \ref{pr:HochSep} and the fact that
$\lambda_\ast$ is a natural map.
\end{proof}

\begin{noname}\label{pa:TauWellDef}
 Our goal now is to construct a cyclic object whose underlying
simplicial object is $\rmd{Z}{\ast}(H,A_B)$. The cyclic operator
can not be taken $t_\ast\tp{H} A_B$, as $t_\ast$ is not in general
a morphism of right $H$--modules. Nevertheless, we shall show that
the linear maps $\tau_\ast:\rmd{Z}{\ast}(H,A_B)\rightarrow
\rmd{Z}{\ast}(H,A_B)$ given by
\begin{equation}\label{eq:tau_n}
\tau_n((h^0,\ldots,h^{n})\otimes_H\overline{a})
 =\tsum(h^{n}a\an{1},h^0,\ldots,h^{n-1})\otimes_H\overline{a}\an{0}.
\end{equation}
define a cyclic structure on $\rmd{Z}{\ast}(H,A_B)$. First let us
construct $\tau_n$, $\forall n\in\mathbb{N}$. Let $\tau'_n$ be the
map $(h^0,\ldots,h^{n})\otimes\overline{a}\mapsto
(h^{n}a\an{1},h^0,\ldots,h^{n-1})\otimes_H\overline{a}\an{0}$. We
 claim that
\[
\tau'_n((h^0,\ldots,h^{n})h\otimes\overline{a})= \tau'_n((h^0,\ldots,h^{n})\otimes
h\overline{a})
\]
for every $h^0,\ldots,h^n,h\in H$ and $\overline{a}\in A_B$. A
simple computation proves us that:
\[
\tau'_n((h^0,\ldots,h^{n})h\otimes\overline{a})=\tsum
(h^{n}h\de{n+1}{a}\an{1},h^0h\de{1},\ldots,h^{n-1}h\de{n})
\otimes_H\overline{a}\an{0}
\]
On the other hand, by (\ref{ec:tau5}), we have $\rho(h\overline{a})=\tsum
h\de{2}\overline{a}\an{0}\otimes h\de{3}{a}\an{1}Sh\de{1}$, so:
\begin{eqnarray*}
&&\tau'_n((h^0,\ldots,h^{n})\otimes h\overline{a})=\tsum
    (h^{n}h\de{3}{a}\an{1}Sh\de{1},h^0,\ldots,h^{n-1})\otimes_H
    h\de{2}\overline{a}\an{0}\\
&&{\hspace*{3cm}=\tsum(h^{n}h\de{3}{a}\an{1}Sh\de{1},h^0,\ldots,h^{n-1})h\de{2}\otimes_H
    \overline{a}\an{0}}\\
&&{\hspace*{3cm}=\tsum(h^{n}h\de{n+3}{a}\an{1}Sh\de{1}h\de{2},h^0h\de{3},
    \ldots,h^{n-1}h\de{n+2})\otimes_H\overline{a}\an{0}}\\
&&{\hspace*{3cm}=\tsum(h^{n}h\de{n+1}{a}\an{1},h^0h\de{1},
    \ldots,h^{n-1}h\de{n})\otimes_H\overline{a}\an{0}}\\
&&{\hspace*{3cm}=\tau'_n((h^0,\ldots,h^{n})h\otimes\overline{a})}.
\end{eqnarray*}
It results that $\tau_n'$ is $H$--balanced, so it induces a map $\tau_n$ verifying
(\ref{eq:tau_n}). Now we can prove that $\rmd{Z}{\ast}(H,A_B)$ is a cyclic object with
respect to $\tau_\ast$.
\end{noname}
\begin{theorem}\label{te:S_Bullet}
$(\rmd{Z}{\ast}(H,A_B),\delta_\ast,\sigma_\ast,\tau_\ast)$ is a cyclic object
isomorphic to $\rmd{Z}{\ast}(A/B)$. In particular, the cyclic homology of
$\rmd{Z}{\ast}(H,A_B)$ is $\rmd{HC}{\ast}(A/B)$.
\end{theorem}

\begin{proof}
We keep the notation from the proof of Proposition \ref{pr:IzoSymp}. It is enough to
prove that
\begin{eqnarray}
    \lambda_\ast t_\ast=\tau_\ast\lambda_\ast\label{eq:taun}.
\end{eqnarray}
because, in this case, $\tau_\ast$ is the unique operator that makes
$(\rmd{Z}{\ast}(H,A_B),\delta_\ast,\sigma_\ast)$ a cyclic object such that
$\lambda_\ast$ becomes an isomorphism of cyclic objects.

Let us denote by $D_l$ (respectively $D_r$) the left (respectively
right) hand side of (\ref{eq:taun}) evaluated at
$a^0\wotB\cdots\wotB a^n$. We get:
\begin{eqnarray*}
D_l&=&\lambda_n(a^n\wotB a^0\cdots\wotB a^{n-1})\\
&=&\tsum(\pr{0}{n-1}{1},\pr{1}{n-1}{2},\ldots,a^{n-1}\an{n},1)
    \otimes_H\overline{a^n\pr{0}{n-1}{0}}\\
&\overset{(\ast)}{=}&\tsum(\pr{0}{n-1}{1},\pr{1}{n-1}{2},\ldots,a^{n-1}\an{n},1)
    \otimes_H\overline{a^n\an{1}(\pr{0}{n}{0}})\\
&=&\tsum(\pr{0}{n-1}{1},\pr{1}{n-1}{2},\ldots,a^{n-1}\an{n},1)a^n\an{1}
    \otimes_H\overline{(\pr{0}{n}{0}})\\
&=&\tsum(\pr{0}{n}{1},\pr{1}{n}{2},\ldots,\pr{n-1}{n}{n},a^n\an{n+1})
    \otimes_H\overline{(\pr{0}{n}{0}}).
\end{eqnarray*}
To deduce the equality denoted by ($\ast$) we used
(\ref{ec:tau4}), and the last equality was obtained by using the
fact that $H^{\otimes (n+1)}$ is a right module via $\Delta_n$.
\begin{eqnarray*}
D_r &=&\tau_n\left(\tsum (\pr{1}{n}{1},\pr{2}{n}{2},\ldots,
    \pr{n-1}{n}{n-1}, a^n\an{n},1)\otimes_H%
    \overline{{a^0}\pr{1}{n}{0}}\right)\\
&=&\tsum (a^0\an{1}\pr{1}{n}{1},\pr{1}{n}{2},\pr{2}{n}{3},\ldots,
    a^{n-1}\an{n}a^{n-1}\an{n}, a^n\an{n+1})\otimes_H%
    \overline{a^0\an{0}\pr{1}{n}{0}}\\
&=&D_l.
\end{eqnarray*}
This sequence of equalities  completes the proof of the theorem.
\end{proof}

\begin{remark}
At a first sight, for proving the previous theorem, we need only that $A_B$ is a left
$H$--module and a right $H$--comodule such that (\ref{ec:tau4}) holds true. In fact,
relation (\ref{ec:tau5}) is required also, because it is equivalent to the fact that
$\tau_1:\rmd{Z}{1}(A_B)\longrightarrow \rmd{Z}{1}(A_B)$ is well--defined.
\end{remark}

\section{Hochschild and cyclic homology of modular crossed modules}
\markboth{\sc Hochschild and cyclic homology of modular crossed
modules}{\sc Hochschild and
cyclic homology of modular crossed modules}%
In this section we shall show that a  cyclic object, similar to
$\rmd{Z}{\ast}(H,A_B)$, can be constructed for each left
$H$--module $M$ on which $H$--coacts in a compatible way. The
particular case $A_B$ suggests the following definition.

\begin{definition}\label{de:CrossedModule}
 Let $H$ be a Hopf algebra and let $M$ be a left $H$--module
which is  a right $H$--comodule with respect to
$\rho:M\longrightarrow M\otimes H$.

a) $M$ is called a \emph{crossed module} if and only if for every $m\in M$ and $h\in H$
\begin{equation}\label{ec:cm}
\rho(hm)=\tsum h\de{2}m\an{0}\otimes h\de{3}m\an{1}Sh\de{1}.
\end{equation}

b) A crossed module $M$ is called \emph{modular}\footnote{The main
results of this paper were presented by the authors in their talks
at the ``First Joint Meeting RSME--AMS", Sevilla, June 2003.
During that meeting Tomasz Brzezi\'{n}ski informed us that modular
crossed modules were also defined in \cite{HKRS}, where they are
called stable anti--Yetter--Drinfeld module.} if and only if for
every $m\in M$
\begin{equation}\label{ec:cmm}
\tsum m\an{1}m\an{0}=m.
\end{equation}

By definition a morphism of crossed modules is a map which is both left $H$--linear and
right $H$--colinear. We obtain a category that will be denoted by $\mathcal{CM}(H)$.
The full subcategory of modular crossed modules will be denoted by $\mathcal{CM}_m(H)$.
\end{definition}

\begin{remark}\label{re:mcm}
If $A/B$ is an $H$--Galois extension then $A_B$ is an example of
modular crossed modules. Other examples are given bellow.
\end{remark}

\begin{example}\label{ex:mcm}
Let $H$ be a Hopf algebra.

a) $H$ is a modular crossed module with respect to the action:
\[
hx=\tsum h\de{2}xSh\de{1},
\]
and the coaction defined by $\Delta$. We shall denote this modular
crossed module by $_{ad}H$. By duality one can construct the
modular crossed module $H^{ad}$. The action on $H^{ad}$ is induced
by the multiplication in $H$, while its coaction is defined by
\[
\rho(h)=\tsum h\de{2}\otimes h\de{3}Sh\de{1}.
\]

b) Connes and Moscovici defined a modular pair to be a pair
$(\sigma,\delta)$  such that $\delta:H\to k$ is a morphism of
algebras, $\sigma\in H$ is a group--like element (that is $\sigma$
is a non--zero element such that
$\Delta(\sigma)=\sigma\tp{}\sigma$) and $\delta(\sigma)=1$.

For such a pair $(\sigma,\delta)$ they constructed a cosimplicial
object $\mathcal H^\natural_{(\sigma,\delta)}$, and they proved
that $\mathcal H^\natural_{(\sigma,\delta)}$ is a cyclic object if
and only if $(L_{\sigma^{-1}}S_\delta)^2=\mathrm{Id}_H$, where
$L_{\sigma^{-1}}$ denotes the left multiplication by
$\sigma^{-1}$, and the twisted antipode $S_\delta:H\to H$ is
defined by $S_\delta(h)=\tsum \delta (h\de{1})S(h\de{2})$. In this
case they call $(\sigma,\delta)$ a modular pair in involution.

One can see easily that $(\sigma,\delta)$ is a modular pair in
involution if and only if the antipode $S$ of $H$ is bijective and
$_\delta k^\sigma\in\mathcal{CM}_m(H^{op\;cop})$. Here
$H^{op\;cop}$ denotes the Hopf algebra $H$ with opposite algebra
and coalgebra structures (and the same antipode $S$), and $_\delta
k^g$ is the one dimensional vector space $k$ on which $H^{op}$
acts via $\delta$ and $H^{cop}$ coacts via $\sigma$.\vsu

 c) We recall that a left Yetter-Drinfeld module over $H$ is
a left $H$--module and a left $H$--comodule $M$ such that
\begin{equation}\label{def:YD}
\rho(hm)=\tsum h\de{1}m\an{-1}Sh\de{3}\tp{}h\de{2}m\an{0}.
\end{equation}
As in the case of crossed modules, we can talk about the category
of Yetter--Drinfeld modules, that will be denoted by $\yd$.

Let us assume that  the antipode of $H$ is involutive, that is
$S^2=\mathrm{Id_H}$. Then we have an isomorphism of categories
$\cm\simeq\yd$. This isomorphism associates to a Yetter--Drinfeld
module $(M,\cdot,\rho)$ the crossed module $(M,\cdot,\rho')$ with
the same left $H$--action and right $H$ coaction given by
$\rho'(m)=\tsum m\an{0}\tp{}Sm\an{-1}$.\vsu

d) Let $f:K\to H$ be a morphism of Hopf algebras and let
$N\in\mathcal{CM}(K)$. Then $\Ind_K^HN:=H\otimes_K N$ is an object
in $\cmm.$ The $H$--module structure on $\Ind_K^HN$ is the
canonical one, while its $H$--comodule is defined by
\[
\rho(h\otimes_K m)=\tsum (h\de{2}\otimes_K m\an{0})\otimes
h\de{3}f(m\an{1})Sh\de{1}.
\]
If $N$ is modular then $\Ind_K^HN$ is modular too.\vsu

e) Let $f:K\to H$ be a morphism of Hopf algebras and let
$M\in\mathcal{CM}(H)$ with comodule structure $\rho_M:M\to
M\otimes H$. Recall that $H$ coacts to the left on $K$ by
$\rho_K:=(f\otimes K)\Delta_K$, and by definition we have
$M\square_H K=\rmd{Ker}{}(\rho_M\otimes K-M\otimes\rho_K$).

Then $\rmd{Res}{K}^H(M):= M\square_H K$ is an object in
$\mathcal{CM}(K)$. The $K$--comodule structure on
$\rmd{Res}{K}^H(M)$ is the canonical one, while its $K$--module
structure is defined by
\[
k(\tsum m_i\square_H x_i)=\tsum f(k\de{2})m_i\square_H
k\de{3}x_iSk\de{1},
\]
where $\tsum m_i\square_H x_i\in M\square_H K$. If $M$ is modular
then $\rmd{Res}{K}^H(M)$ is modular too. Note that
$\rmd{Res}{K}^H(-)$ is a right adjoint functor of
$\Ind_K^H(-)$.\vsu

f) Let us consider the case when $H=kG$. A $kG$--comodule
$(M,\rho)$ is a vector space together with a decomposition
$M=\bigoplus_{x\in{G}}M_x$. The subspace $M_x$ contains all
elements $m\in M$ such that $\rho(m)=m\tp{}x$. Hence a
$kG$--crossed module is a $kG$--module $M$ having a decomposition
$M=\bigoplus_{x\in{G}}M_x$ such that
\[
 g{M_x}\subseteq{M_{gxg^{-1}}},\quad\forall g,x\in{G}.
\]

For $g\in{G}$ let $[g]$ denote the conjugacy class of $g$. Let
$T(G)$ be the set of all conjugacy classes in $G$ and let $t(G)$
be a coset for $T(G)$, so every element in $T(G)$ is the conjugacy
class of one and only one element in $t(G)$.

Suppose that $M$ is a crossed module. If  $g\in{t(G)}$ we define
$M_{[g]}:=\bigoplus_{x\in[g]}M_x$. Obviously $M_{[g]}$ is a
crossed submodule of $M$ as $x{M_{[g]}}\subseteq{M_{[g]}}$ for
every $x\in{G}$.
\end{example}

\begin{proposition}\label{pr:Crossed_G}
Let $M\in\mathcal{CM}(kG)$. Then there is a canonical isomorphism of crossed modules:
\begin{equation}\label{eq:Crossed_G}
 M\simeq\oplus_{g\in{t(G)}}M_{[g]}
\end{equation}
and $M_{[g]}\simeq\rmd{Ind}{kG_g}^{kG}M_g$, where $G_g$ is the
centralizer of $g$ in $G$. $M$ is modular if and only if $gm=m$
for every $g\in G$ and $m\in M_g$.
\end{proposition}
\begin{proof}

Obviously we have the decomposition (\ref{eq:Crossed_G}) since $M$
is a $kG$--comodule. We have to show that
$M_{[g]}\simeq\rmd{Ind}{kG_g}^{kG}M_g$ for every $g\in{t(G)}$.

Let $g\in t(G)$. First let us remark that $x{M_g}=M_g$ for every
$x\in{}G_g$, so $G_g$ acts on $M_g$. If
$\varphi:kG\otimes{M_g}\longrightarrow{M_{[g]}}$ denotes the map
induced by the module structure of $M$, then $\varphi$ is
surjective since the multiplication by $x\in{G}$ is a bijective
$k$--linear map from $M_g$ to $M_{xgx^{-1}}$. The map $\varphi$
induces a surjective morphism of crossed modules
$\overline{\varphi}:\rmd{Ind}{kG_g}^{kG}M_g\to{M_{[g]}}$. Let $R$
be a coset for $(G/G_g)_r=\{xG_g\mid\;x\in{G}\}$. The set
$\{e_x\mid\;x\in{R}\}$ is a basis of $kG$ as a right
$kG_g$--module, where $e_x=\sum_{y\in{}G_g} xy$. Hence an element
$u\in{kG}\otimes_{kG_g}M_g$ can be written uniquely as a sum
\[
 u=\tsum_{x\in{R}}e_x\otimes_{kG_g}m_x,
\]
with $m_x\in{M_g}$. If $\overline{\varphi}(u)=0$, then
$\sum_{x\in{R}}\sum_{h\in{}G_g}xh{m_x}=0$. But $h{m_x}\in{M_g}$ as
$h\in{}G_g$. On the other hand, for $x\neq{y}$ in $R$ we have
$xgx^{-1}\neq{ygy^{-1}}$. Thus $\sum_{h\in{}G_g}h{m_x}=0$ for all
$x\in{R}$. We deduce that
\[
 u=\tsum_{x\in{R}}\tsum_{h\in{}G_g}xh\otimes{m_x}
 =
 \tsum_{x\in{R}}x\otimes\tsum_{h\in{}G_g}hm_x=0,
\]
therefore $\overline{\varphi}$ is injective too.
\end{proof}

\begin{remark}
Let $\mathcal{C}:=\prod_{x\in t(G)}{}_{G_x}\mathcal{M}$. The
previous proposition shows us that $\mathcal{CM}(kG)$ is
equivalent to $\mathcal{C}$.
\end{remark}

We are going to investigate some other homological properties of
$\mathcal{CM}(H)$ and $\mathcal{CM}_m(H)$.

\begin{proposition}
$\mathcal{CM}(H)$ is an abelian category and $\mathcal{CM}_m(H)$ is a closed class in
$\mathcal{CM}(H)$ with respect to quotients, subobjects and direct sums. In particular,
$\mathcal{CM}_m(H)$ is abelian too.
\end{proposition}
\begin{proof}
Straightforward, left to the reader.
\end{proof}

\begin{lemma}
Let $M$ be a crossed module over $H$ and let $u_M:M\to M$ be the
map
\[
u_M(m)=\tsum{m}\an{1}{m}\an{0}.
\]
Then $u_M$ is a morphism of crossed modules. The kernel and the
cokernel of $u_M-\rmd{Id}{M}$ are modular crossed modules.
\end{lemma}

\begin{proof}
By relation (\ref{ec:cm}) it follows immediately that $u_M$ is
left $H$--linear. Since
\begin{equation}\label{eq:subcomodule}
 \rho(\tsum m\an{1}m\an{0})
 =\tsum{m}\de{3}{m}\an{0}\otimes{m}\an{4}m\an{1}Sm\an{2}
 =\tsum{m}\an{1}{m}\an{0}\otimes{m}\an{2}
\end{equation}
it results that $u_M$ is $H$--colinear too. Obviously the kernel
and the cokernel of $u_M-\rmd{Id}{M}$ are modular.
\end{proof}

\begin{proposition}
Let $M$ be a modular crossed module over a Hopf algebra $H$. The
forgetful functor $_HF:\cmm\longrightarrow{_H\mm}$ has a right
adjoint $_HG$.
\end{proposition}
\begin{proof}
Let $F':\mm^H\longrightarrow\mm_k$ be the functor that forgets the
comodule structure. It is well known that
$G':\mm_k\longrightarrow\mm^H$, $G'(V)=V\otimes{H}$ is a right
adjoint functor of $F'$, where the comodule structure of $G'(V)$
is defined by $V\otimes\Delta_H$. The adjunction is given by the
isomorphisms:
\begin{equation}\label{eq:alphabeta}
\begin{array}{ll}
 \alpha_{M,V}:\Hom_k(F'(M),V)\longrightarrow\Hom^H(M,G'(V)),
 \quad
 \alpha_{M,N}(f)=(f\otimes{H})\rho_M,\\
 \beta_{M,V}:\Hom^H(M,G'(V))\longrightarrow\Hom_k(F'(M),V),
 \quad
 \beta_{M,V}(g)=(V\otimes\varepsilon)g.
\end{array}
\end{equation}
If we take $N$ to be a left $H$--module, then $G'(N)$ becomes a
left $H$--module with the structure
\[
 h({n}\otimes{x}):=\tsum h\de{2}n\otimes{h}\de{3}xSh\de{1},
\]
so that $G'(N)$ is a crossed module over $H$. Therefore $G'$
induces a functor, denoted also by $G'$, from $_H\mm$ to $\cm$. In
general $G'(N)$ is not modular, but by the previous Lemma
\[
 {}_HG(N):=\rmd{Ker}{}(u_{G'(N)}-\rmd{Id}{G'(N)})
\]
is a modular crossed module. The correspondence
$N\longmapsto{_HG(N)}$ defines the functor that we are looking
for. The adjunction $_HF\adjoint{_HG}$ is defined by the
restrictions of $\alpha$ and $\beta$ defined in
(\ref{eq:alphabeta}).
\end{proof}

\begin{remark}
The unit of the adjunction $\sigma_M:M\longrightarrow{_HG(M)}$ is
the corestriction of $\rho_M$, the map defining the comodule
structure, to $_HG(M)\subseteq{M\otimes{H}}$.
\end{remark}

Note that $\sigma_M$ is a monomorphism, because $\rho_M$ is so,
thus by embedding $M$ into an injective module $I$ it results that
any modular crossed module can be embedded into an injective one,
namely $_HG(I)$. For future references we state this result in the
following proposition.

\begin{proposition}
The category $\cmm$ is an abelian category with enough injective
objects.
\end{proposition}
\begin{noname}
In a similar way one can construct a left adjoint $G^H$ of
$F^H:\cmm\to\calM^H$, the functor that forgets the module
structure. For every right $H$--comodule $N$ we can regard
$H\otimes N$ as a crossed module with respect to the structures
\begin{eqnarray*}
&&h(x\otimes m)=(hx)\otimes m,\\
&&\rho(x\otimes m)=\tsum (x\de{2}\otimes m\an{0})\otimes
x\de{3}m\an{1}Sx\de{1}.
\end{eqnarray*}
For every right comodule $N$ we set
$G^H(N)=\rmd{Coker}{}(u_{H\otimes N}-\rmd{Id}{H\otimes N})$. By
construction $G^H(N)$ is a modular crossed module, and the
correspondence $N\longmapsto G^H(N)$ defines a functor, which is
left adjoint of $F^H$. The details of the proof are left to the
reader.
\end{noname}

\begin{remark}
The functors $_HG$ and $G^H$ provides new examples of crossed
modules. If $_\varepsilon k$  denotes the trivial $H$--action on
$k$ then $_HG(_\varepsilon k)={}_{ad}H$. For the definition of
$_{ad}H$ see Example \ref{ex:mcm}(a).
\end{remark}

Now we are going to associate to every modular crossed module $M$
a new cyclic object. Since $M$ is in particular a left $H$--module
one can consider the simplicial object
$\rmd{Z}{\ast}(H,M)=H^{\otimes(\ast+1)}\otimes_H M$ that was
constructed in Corollary \ref{co:S(M)}. For reader convenience we
recall the definition of $(\delta_i)_{0\leq i\leq n}$ and
 $(\sigma_i)_{0\leq i\leq n}$, respectively the face and degeneracy
maps of $\rmd{Z}{\ast}(H,M)$ in degree $n$.
\begin{eqnarray}
 \delta_i((h^0,\ldots,h^{n})\otimes_Hm)
 &=&\varepsilon(h^{i})(h^0,\ldots,\widehat{h^{i}},\ldots,h^{n})
 \otimes_Hm,\label{eq:delta}\\
 \sigma_i((h^0,\ldots,h^{n})\otimes_Hm)
 &=&(h^0,\ldots,\Delta{h^{i}},\ldots,h^{n})\otimes_Hm.\label{eq:sigma}
\end{eqnarray}
To define the cyclic operator
$\tau_n:\rmd{Z}{n}(H,M)\longrightarrow\rmd{Z}{n}(H,M)$ we use the
comodule structure of $M$.
\begin{equation}\label{eq:tau}
\tau_n((h^0,\ldots,h^{n})\otimes_H{m})
 =\tsum(h^{n}m\an{1},h^0,\ldots,h^{n-1})\otimes_H{m}\an{0}.
\end{equation}
Using the compatibility relation from the definition of crossed
modules  we can prove that $\tau_\ast$ is well--defined as in
(\ref{pa:TauWellDef}). Note that the cyclic operator
(\ref{eq:tau_n}) is a particular case of (\ref{eq:tau}).
\begin{theorem}\label{te:CyclicObject}
Let $M$ be a modular crossed module over a Hopf algebra $H$. Then
$(\rmd{Z}{\ast}(H,M),\delta_\ast,\sigma_\ast,\tau_\ast)$ is a
cyclic object.
\end{theorem}
\begin{proof}
We know from Corollary \ref{co:S(M)} that $\rmd{Z}{\ast}(H,M)$ is
a simplicial object, and we have seen that the cyclic operator is
well--defined. Therefore we have to check the following relations:
\[
\delta_i\tau_n=\tau_{n-1}\delta_{i-1},\qquad
\sigma_i\tau_n=\tau_{n+1}\sigma_{i-1},\qquad
\delta_0\tau_n=\delta_n,\qquad
\sigma_0\tau_n=\tau_{n+1}^2\sigma_n.
\]
where $0<i\leq n$. Their proofs are consequences of the
computations bellow.

\begin{eqnarray*}
&&\delta_i\tau_n((h^0,\ldots,h^{n})\tp{H}m)=\delta_i(\tsum
    (h^{n}m\an{1},h^0,\ldots,h^{n-1})\tp{H}m\an{0})\\
&&\hspace*{2cm}=\tsum\varepsilon(h^{i-1})(h^{n}m\an{1},h^0,\ldots,
    \widehat{h^{i-1}},\ldots,h^{n-1})\tp{H}m\an{0},\\
&&\tau_{n-1}\delta_{i-1}((h^0,\ldots,h^{n})\tp{H}m)
    =\tsum\tau_{n-1}(\varepsilon(h^{i-1})(h^0,\ldots,
    \widehat{h^{i-1}},\ldots,h^{n})\tp{H}m)\\
&&\hspace*{2cm}=\tsum
\varepsilon(h^{i-1})(h^{n}m\an{1},h^0,\ldots,
    \widehat{h^{i-1}},\ldots,h^{n-1})\tp{H}m\an{0}.
\end{eqnarray*}
\begin{eqnarray*}
&&\sigma_i\tau_n((h^0,\ldots,h^{n})\tp{H}m)=\sigma_i (\tsum
    (h^{n}m\an{1},h^0,\ldots,h^{n-1})\tp{H}m\an{0})\\
&&\hspace*{2cm}=\tsum
    (h^{n}m\an{1},h^0,\ldots,\Delta(h^{i-1}),\ldots,h^{n-1})\tp{H}m\an{0},\\
&&\tau_{n+1}\sigma_{i-1}((h^0,\ldots,h^{n})\tp{H}m)=
    \tau_{n+1}((h^0,\ldots,\Delta(h^{i-1}),\ldots,h^{n})\tp{H}m)\\
&&\hspace*{2cm}=\tsum
    (h^{n}m\an{1},h^0,\ldots,\Delta(h^{i-1}),\ldots,h^{n-1})\tp{H}m\an{0}.
\end{eqnarray*}
\begin{eqnarray*}
&& \delta_0\tau_n((h^0,\ldots,h^{n})\tp{H}m)=\delta_0(\tsum
    (h^{n}m\an{1},h^0,\ldots,h^{n-1})\tp{H}m\an{0})\\
&&\hspace*{2cm}=\tsum\varepsilon(h^{n}m\an{1})(h^0,\ldots,h^{n-1})\tp{H}m\an{0}\\
&&\hspace*{2cm}=\tsum\varepsilon(h^{n})(h^0,\ldots,h^{n-1})\tp{H}m\\
&&\hspace*{2cm}=\delta_n((h^0,\ldots,h^{n})\tp{H}m).
\end{eqnarray*}
\begin{eqnarray*}
&&\sigma_0\tau_n((h^0,\ldots,h^{n})\tp{H}m)=\sigma_0(\tsum
    (h^{n}m\an{1},h^0,\ldots,h^{n-1})\tp{H}m\an{0})\\
&&\hspace*{2cm}=\tsum (h^{n}\de{1}m\an{1},h^{n}
    \an{2}m\an{2},h^0,\ldots,h^{n-1})\tp{H}m\an{0},\\
&&\tau_{n+1}^2\sigma_n((h^0,\ldots,h^{n})\tp{H}m)=\tau_{n+1}^2(\tsum(
    h^0,\ldots,h^{n-1},h^{n}\de{1},h^{n}\de{2})\tp{H}m)\\
&&\hspace*{2cm}=\tau_{n+1}(\tsum(h^{n}\de{2}m\an{1},
    h^0,\ldots,h^{n-1},h^{n}\de{1})\tp{H}m\an{0})\\
&&\hspace*{2cm}=\tsum (h^{n}\de{1}m\an{1},h^{n}
    \de{2}m\an{2},h^0,\ldots,h^{n-1})\tp{H}m\an{0}.
\end{eqnarray*}
It remains to prove that $\tau_n^{n+1}=\rmd{Id}{}$. But
\begin{equation*}
\tau_{n}^{n+1}((h^0,\ldots,h^{n})\tp{H}m)=
    \tsum(h^{0}m\an{1},\ldots,h^{n}m\an{n+1})\tp{H}m\an{0}
\end{equation*}
equation which implies
\begin{equation*}
\tau_{n}^{n+1}((h^0,\ldots,h^{n})\tp{H}m)
=\tsum(h^{0},\ldots,h^{n})\tp{H}m\an{1}m\an{0}
=\tsum(h^{0},\ldots,h^{n})\tp{H}m
\end{equation*}
Hence $\tau_n^{n+1}=\rmd{Id}{}$.
\end{proof}

\begin{remark}
By definition, an $H$--module coalgebra is a coalgebra $C$ endowed
with a right $H$--module structure such that both the
comultiplication and the counit of $C$ are morphisms of
$H$--modules ($H$ acts diagonally on $C\otimes C$ and trivially on
$k$).

For every $H$--module coalgebra $C$ and every $M\in\cmm$ we can
modify the definition of
$(\rmd{Z}{\ast}(H,M),\delta_\ast,\sigma_\ast,\tau_\ast)$ to obtain
a new cyclic object as follows. First, let
\[
\rmd{Z}{\ast}(C,M):=C^{\otimes(\ast+1)}\otimes_H M
\]
and then define $\delta_\ast,\sigma_\ast,\tau_\ast$ as in
(\ref{eq:delta}), (\ref{eq:sigma}), (\ref{eq:tau}), of course,
using the comultiplication and the counit of $C$ this time. Then
it is easy to see that
$(\rmd{Z}{\ast}(C,M),\delta_\ast,\sigma_\ast,\tau_\ast)$ is a
cyclic object.
\end{remark}

\begin{definition}
The complex associated to the simplicial object
$\rmd{Z}{\ast}(H,M)$ will be denoted by $\rmd{C}{\ast}(H,M)$. For
the Tsygan double complex, the Hochschild homology and the cyclic
homology of $\rmd{Z}{\ast}(H,M)$ we will use the notation
$\rmd{CC}{\ast\ast}(H,M)$, $\rmd{HH}{\ast}(H,M)$ and
$\rmd{HC}{\ast}(H,M)$, respectively.
\end{definition}

\begin{proposition}
The functor that associates to each $M\in\mathcal{CM}_m(H)$ the
cyclic object $\rmd{Z}{\ast}(H,M)$ is an exact functor that
commutes with direct sums.
\end{proposition}
\begin{proof}
By definition we have
$\rmd{Z}{\ast}(H,M)=\wrmd{Z}{\ast}(H)\otimes_HM$ as vector spaces.
The proposition follows by the fact that $\wrmd{Z}{\ast}(H)$ is a
free right $H$--module with respect to the diagonal action. The
other assertion is trivial since the tensor product commutes with
direct sums.
\end{proof}

$H$-module algebras are defined by duality from $H$--comodule
algebras. For an arbitrary  modular pair in involution
$(\delta,\sigma)$ and an arbitrary $H$--module algebra $A$, Connes
and Moscovici defined a $\delta$--invariant $\sigma$--trace to be
a morphism of $H$--modules $\tau:A\to {}_\delta k$ such that
$\tau(ab)=\tau(b\sigma(a))$, for all $a,b\in A$. They proved that
such a $\tau$ defines a morphism of cocyclic objects from
$H^\#_{(\delta,\sigma)}$ to $A^\#$, the usual cocyclic module
associated to the algebra $A$. Hence, for $\tau$ as above, there
is a morphism
\[
\gamma_\tau^\ast:\mathrm{HC}^{\ast}_{(\delta,\sigma)}(H)\to\rmu{HC}{\ast}(A).
\]
For details the reader is referred to \cite{CM1,CM2}. Now we want
to show that a similar construction exists in our setting. First
let us define the appropriate type of trace that we shall use.

\begin{definition}\label{de:trace}
Let $A$ be an $H$-comodule algebra and let $M$ be a modular
crossed module. An $H$--coinvariant $M$--trace of $A$ is a
morphism  $tr_M:A\to M$ of $H$--comodules  such that
\[
tr_M(ax)=\tsum a\an{1}tr_M(xa\an{0}).
\]
\end{definition}

\begin{noname}
The map $\beta_n(A):A^{\otimes_B(n+1)}\to A\otimes H^{\otimes n}$
defined in (\ref{nn:Galois}) exists for an arbitrary $H$--comodule
algebra $A$, but in general it is not bijective. Let us consider
the following sequence of $k$--linear maps:
\[
A^{\otimes(n+1)}\stackrel{can}{\longrightarrow}A^{\otimes_B(n+1)}
\stackrel{\beta_n(A)}{\longrightarrow}A\otimes H^{\otimes
n}\stackrel{\tau\otimes\rmd{Id}{}}{\longrightarrow}M\otimes
H^{\otimes n}\simeq H^{\otimes n}\otimes M\simeq
H^{\otimes(n+1)}\otimes_H M,\
\]
where the first isomorphism is the canonical flip, while the
second one is $\varphi_n\otimes_H\rmd{Id}{M}$, see Remark
\ref{le:030115} for the definition of $\varphi_n$. The composition
of these morphisms, denoted $\gamma_n^{tr_M}$, maps
$x=a^0\otimes\cdots\otimes a^n$ to
\begin{equation}
\gamma_n^{tr_M}(x)=\tsum (\pr{1}{n}{1},\pr{2}{n}{2},\ldots,
    a^{n-1}\an{n-1}a^{n}\an{n-1}, a^n\an{n},1)\otimes_H%
    tr_M({a^0}\pr{1}{n}{0}).
\end{equation}
\end{noname}

\begin{theorem}
Let $A$ be an $H$--comodule algebra and let $M$ be a modular
crossed module over $H$. If $tr_M:A\to M$ is an $H$--coinvariant
$M$--trace then
\[
\gamma_\ast^{tr_M}:\rmd{Z}{\ast}(A/k)\to\rmd{Z}{\ast}(H,M)
\]
is a morphism of cyclic objects that induces a map, denoted also
by $\gamma_\ast^{tr_M}$,
\[
\gamma_\ast^{tr_M}:\rmd{HC}{\ast}(A)\to\rmd{HC}{\ast}(H,M).
\]
\end{theorem}
\begin{proof}
One can proceed as in the proof of the fact that
$\lambda_\ast:\rmd{Z}{\ast}(A/B)\to\rmd{Z}{\ast}(H,A_B)$ is an
isomorphism of cyclic objects, see the proofs of Proposition
\ref{pr:IzoSymp} and Theorem \ref{te:S_Bullet}. The details are
left to the reader.
\end{proof}

\section{Cyclic homology of induced modular crossed modules}
\markboth{\sc Cyclic homology of induced modular crossed
modules}{\sc Cyclic homology of induced modular crossed modules}

Throughout this section we whall assume that the characteristic of
$k$ is zero. Let $H$ be a Hopf algebra over $k$, and let $K$ be a
Hopf subalgebra of $H$. Let $M$ be a given modular crossed
$K$--module $M$.

The main goal of this section is to compute, under some
assumptions on $K$ and $M$, the cyclic homology of $H$ with
coefficients in $\rmd{Ind}{K}^HM$. As applications we shall
compute the cyclic homology of strongly graded algebras, group
algebras and quantum tori. We start by proving a variant of
Shapiro's Lemma for cyclic homology.

\begin{proposition}\label{pr:229}
Let $H$ be a Hopf algebra and let $K$ be a Hopf subalgebra of $H$.
If $M$ is a modular crossed $K$--module  and $H$ is a flat right
$K$--module, then
\[
 \rmd{HH}{\ast}(H,\Ind_K^HM)\simeq\rmd{HH}{\ast}(K,M)
 \mbox{ and }
 \rmd{HC}{\ast}(H,\Ind_K^HM)\simeq\rmd{HC}{\ast}(K,M).
\]
\end{proposition}
\begin{proof}
The cyclic object associated to $\Ind_K^HM$ can be identified as follows:
\[
 \rmd{Z}{\ast}(H,\Ind_K^HM)\simeq{H^{\otimes{\ast+1}}}\otimes_KM,\
\]
where the formulas for the face, degeneracy and cyclic operators
are obtained from (\ref{eq:delta}), (\ref{eq:sigma}) and
(\ref{eq:tau}) by replacing $\otimes_H$ with $\otimes_K$
everywhere.

Since $H$ is a flat $K$--module, it results that
$\wrmd{C}{\ast}(H)\stackrel{\varepsilon}{\longrightarrow}{k}$ is
also a flat resolution of $k$ over $K$, therefore
\[
 \rmd{HH}{n}(H,\Ind_K^HM)\simeq\Tor{K}{n}(k,M)\simeq\rmd{HH}{n}(K,M).
\]
The isomorphism above can be described explicitly as follows. The
inclusion map $i:K\longrightarrow{H}$ induces a morphism
$i_\ast:\rmd{\widetilde{C}}{\ast}(K)\longrightarrow
\rmd{\widetilde{C}}{\ast}(H)$ of resolutions of $k$ over $K$.
Hence $i_\ast\otimes_KM$ is a quasi--isomorphism of complexes. Of
course $i_\ast\otimes_KM$ is a morphism of cyclic objects, so by
the 5--Lemma we get:
\[
 \rmd{HC}{\ast}(K,M){\simeq}\rmd{HC}{\ast}(H,\Ind_K^HM).
\]
Note that the last isomorphism is also induced by
$i_\ast\otimes_KM$.
\end{proof}

In the particular case when $H$ is cocommutative we can compute
easily the cyclic homology of a crossed module $M$ with trivial
coaction, i.e. $\rho(m)=m\otimes1$, for all $m\in{M}$.

\begin{theorem}\label{th:228}
Let $H$ be a cocommutative Hopf algebra. If $M$ is a modular
crossed module over $H$ with trivial coaction, then
\[\textstyle
 \rmd{HC}{\ast}(H,M)=\bigoplus_{i\geq0}\Tor{H}{\ast-2i}(k,M).
\]
\end{theorem}
\begin{proof}
Let $\wrmd{Z}{\ast}(H)$ be the cyclic object constructed in
Proposition~\ref{pr:C_*(H)}. Since $H$ is cocommutative $t_\ast$
is right $H$--linear, hence it makes sense $t_\ast\otimes_HM$.
Because $H$ coacts trivially on $M$ then
$\tau_\ast=t_\ast\otimes_HM$, where $\tau_\ast$ is the cyclic
operator of $\rmd{Z}{\ast}(H,M)$. Thus:
\[
 \rmd{Tot}{}\rmd{CC}{\ast\ast}(H,M)\simeq\rmd{Tot}{}\wccd(H)\otimes_HM,
\]
On the other hand $\wrmd{CC}{\ast\ast}(H)$ is a resolution (in the
sense of hyperhomolgy) of
\[
\rmd{C}{\ast}:\qquad 0\longleftarrow k\longleftarrow
0\longleftarrow k\longleftarrow\cdots
\]
($C_n=0$ if either $n<0$ or $n$ is odd). As in the proof of
Karoubi's theorem, see for example \cite[Theorem 9.7.1]{We}, it
follows that $\rmd{HC}{\ast}(H,M)$ is the hypertor
$\mathbb{T}\rmd{or}{\ast}^H(\rmd{C}{\ast},M)$. It is easy to see
that another resolution of $\rmd{C}{\ast}$ is
\[
 0\longleftarrow \wrmd{C}{\ast}(H)\longleftarrow 0\longleftarrow
\wrmd{C}{\ast}(H)\longleftarrow\cdots.
\]
If we use this resolution to compute
$\mathbb{T}\rmd{or}{\ast}^H(\rmd{C}{\ast},M)$ we immediately get
that
\[\textstyle
\mathbb{T}\rmd{or}{\ast}^H(\rmd{C}{\ast},M)=\bigoplus_{i\geq0}\Tor{H}{\ast-2i}(k,M),
\]
equality that proves the theorem.
\end{proof}

\begin{noname}\label{nn:GalExt}
A Hopf subalgebra $K$ of $H$ is called normal if $K_+H$ is a
two--sided ideal in $H$, where $K_+$ denoted the kernel of
$\varepsilon_K$. In fact $K_+H$ is a Hopf ideal in this case, so
$\overline{H}:=H/K_+H$ is a Hopf algebra called the \emph{quotient
Hopf algebra} of $H$ by $K$. It is well-known that, whenever $H$
is a faithfully flat left (or right) $K$--module, then
$K\subseteq{H}$ is an $\overline{H}$--Galois extension \cite[p.
197]{Sch}. Of course the coaction of $\overline{H}$ on $H$ is
defined by the morphism of Hopf algebras
$\pi:H\longrightarrow\overline{H}$.
\end{noname}

\begin{lemma}\label{le:GalExt}
Let $H$ be a Hopf algebra with bijective antipode and let $K$ be a
semisimple normal Hopf subalgebra of $H$. Then $K\subseteq H$ is
an $\overline{H}$--Galois extension.
\end{lemma}

\begin{proof}
Each semisimple Hopf algebra is separable \cite{St}. In
particular, the antipode of $K$ is bijective, as any separable
algebra (over a field) is finite dimensional. Now we can apply
\cite[Corollary 2.9]{Ma} to show that $H$ is a faithfully flat
$K$--module, so the lemma follows by (\ref{nn:GalExt}).
\end{proof}

\begin{noname}\label{nn:Ulbrich}
Let us suppose that $K$ is a normal Hopf subalgebra of $H$ such
that the extension $K\subseteq H$ is $\overline{H}$--Galois. Let
$M$ be an $H$--crossed module. We shall regard $M$ as an
$(H,H)$--bimodule with trivial right  action (via the counit
$\varepsilon$ of $H$). Therefore $(\rmd{Id}{H},\pi,\rmd{Id}{M})$
is a morphism in $\mathfrak C$ from $(H,k,H,M)$ to
$(H,K,\overline{H},M)$, where $\pi$ is the canonical projection
(see Remark \ref{re:Cat_C} for the definition of $\mathfrak C$).
Obviously  $\overline{M}:=M/K_+M$ is a modular crossed module over
$\overline{H}$ with respect to the quotient structures. As ${M}$
is a bimodule with trivial right $\overline{H}$ action it follows
that $M_K:=M/[M,K]$ is equal to $\overline{M}$ as a vector space.
In fact we have more than that. Since by assumption $H/K$ is
$\overline{H}$--Galois it makes sense to speak about the
Ulbrich--Miyashita action of $\overline{H}$ on ${M}_{K}$. As the
canonical map $\beta:H\otimes_K H\to H\otimes \overline{H}$ from
the definition of Galois extensions has the property
\[
\beta(\tsum Sh\de{1}\otimes_K h\de{2})=1\otimes \overline{h}
\]
it results easily that $M_K$ (with Ulbrich-Miyashita action) can
be identified as an $\overline{H}$--module with $\overline{M}$.

Now we can prove the following result.
\end{noname}

\begin{proposition}\label{pr:K=sep}
Let $H$ be a Hopf algebra with bijective antipode. Suppose that
$K$ is a normal semisimple Hopf subalgebra of $H$. For $M\in\cmm$
let $\overline{M}$ denote $M/K_+M$. Then
\[
\begin{array}{ll}
 \rmd{HH}{\ast}(H,M)\simeq\rmd{HH}{\ast}(\overline{H},\overline{M})
 \quad\mbox{ and }
 \quad
 \rmd{HC}{\ast}(H,M)\simeq\rmd{HC}{\ast}(\overline{H},\overline{M}).
\end{array}
\]
\end{proposition}

\begin{proof}
By Lemma \ref{le:GalExt} the extension $K\subseteq H$ is
$\overline{H}$--Galois, and we have already noticed in
(\ref{nn:GalExt}) that a semisimple Hopf algebra is separable. By
the functorial character (both in $H$ and $M$) of the cyclic
object $\mathrm{Z}_\ast(H,M)$ there is a canonical morphism of
cyclic objects
\[
\varphi_\ast:\rmd{Z}{\ast}(H,M)\longrightarrow\rmd{Z}{\ast}(\overline{H},M_K),
\]
that is induced by the projections $H\to\overline{H}$ and $M\to
M_K$. Obviously $\varphi_\ast$ is a morphism of complexes from
$\rmd{C}{\ast}(H,M)$ to $\rmd{C}{\ast}(\overline{H},M_K)$. By
Corollary~\ref{co:LambdaFunct}, it results that $\varphi_\ast$ is
a quasi--isomorphism. We have already proved in (\ref{nn:Ulbrich})
that $M_K\simeq\overline{M}$ as $\overline{H}$--modules. Hence
\[
 \rmd{HH}{\ast}(H,M)\simeq\rmd{HH}{\ast}(\overline{H},\overline{M}).
\]
The isomorphism for cyclic homology follows by using the
SBI--sequence and 5--Lemma.
\end{proof}

\begin{corollary}\label{co:232}
Keep the notation and assumptions from the preceding proposition.
If in addition $\overline{H}$ is cocommutative and
$\rho(M)\subseteq{M\otimes{K}}$ then
\[\textstyle
 \rmd{HC}{\ast}(H,M)\simeq\bigoplus_{i\geq 0}\Tor{\overline{H}}{*-2i}(k,\overline{M}).
\]
\end{corollary}
\begin{proof}
Since the $\overline{H}$--coaction on $\overline{M}$ is trivial
and $\overline{H}$ is cocommutative, we apply
Proposition~\ref{pr:K=sep} and Theorem~\ref{th:228}.
\end{proof}

If $K$ is not separable, then the cyclic cohomology of a modular
crossed module can not be computed in a similar way.

Recall that a Hopf algebra $H$ is by definition pointed if it has
only one dimensional subcoalgebras. Since there is a one--to--one
correspondence between simple right $H$--comodules and simple
subcoalgebras it follows that for each simple comodule $M$ there
is a group--like elements $g\in H$ such that
\[
\rho(m)=m\otimes g,\quad\forall m\in M.
\]

\begin{proposition}\label{pr:030407}
Let $H$ be a cocommutative pointed Hopf algebra. Suppose that
$x\in{H}$ is a central group--like element of infinite order. If
$M\in \cmm$ and $\rho(m)=m\otimes{x}$, for all $m\in{M}$, then
\[
 \rmd{HC}{\ast}(H,M)\simeq\Tor{\overline{H}}{\ast}(k,M),
\]
where $\overline{H}:=H/(1-x)H$.
\end{proposition}
\begin{proof}
Let us remark first that $x{m}=m$, for every $m\in{M}$, since $M$
is modular. Hence $M$ is a $\overline{H}$--module. Now let us
define $t_n:H^{\otimes{(n+1)}}\longrightarrow{H^{\otimes{(n+1)}}}$
by
\[
 t_n(h^0,\ldots,h^{n})=(h^{n}x,h^0,\ldots,h^{n-1}).
\]
Since $H$ is cocommutative, $t_n$ is right $H$--linear (with
respect to the diagonal action) and the cyclic operator
$\tau_\ast$ of $\mathrm{Z}_{\ast}(H,M)$ satisfies the relation
\[
 \tau_\ast=t_\ast\otimes_HM.
\]
Trivially $t_n^{n+1}(u)=ux$, for all $u\in{H^{\otimes{n+1}}}$, and
$t_n$ verifies all other properties of cyclic operators. Moreover,
by Proposition~\ref{pr:C_*(H)}(b), the complex
$\widetilde{\rmd{C}{\ast}}(H)$ is a free resolution of $k$, and
$H$ is free over $k[x]$, see  \cite[p.271]{Ra}. As the order of
$x$ is infinite it results that $k[x]$ is a domain, so $1-x$ is
not a zero--divisor in $H$. Now we can conclude by applying the
following lemma.
\end{proof}

\begin{lemma}
Let $A$ be an algebra over a field $k$ of characteristic $0$ and
let $M$ be a left $A$--module. Suppose that $a\in{A}$ is a central
element such that $1-a$ is not a zero--divisor. Denote $A/(1-a)A$
by $\overline{A}$. Let $(\mathrm{Z}_{\ast},\partial_\ast,s_\ast)$
be a simplicial object in the category of right $A$--modules such
that each $\mathrm{Z}_n$ is free over $A$. Assume that
$t_\ast:Z_\ast\to Z_\ast$ is a morphism of $A$--modules such that
$t^{\ast+1}_\ast=\mathrm{Id}_{Z_\ast}{a}$ and
$(\mathrm{Z}_{\ast},\partial_\ast,s_\ast,t_\ast)$ satisfies all
other properties of cyclic objects. If the complex
$(\mathrm{Z}_{\ast},d_\ast)$ associated to
$(\mathrm{Z}_\ast,\partial_\ast,s_\ast)$ is acyclic and
$(1-a)M=0$, then $\mathrm{Z}_\ast\otimes_A{M}$ is a cyclic object
and
\[
 \rmd{HC}{\ast}(\mathrm{Z}_{\ast}\otimes_A{M})
 \simeq\Tor{\overline{A}}{\ast}(\rmd{H}{0}(Z_\ast),{M}).
\]

\end{lemma}
\begin{proof}
Since, by assumption, we have
$t^{\ast+1}_\ast=\mathrm{Id}_{Z_\ast}{a}$ and $(1-a)M=0$ then
obviously $(\mathrm{Z}_{\ast}\otimes_A {M},\partial_\ast\otimes_A
{M},s_\ast\otimes_A {M},t_\ast\otimes_A {M})$ is a cyclic object.
Let $(\mathrm{Z}'_\ast,d'_\ast)$ be the complex defined by
$Z'_\ast=Z_\ast$  and
$d'_\ast=\sum_{i=0}^{\ast-1}(-1)^{i}\partial_i$. Its homology is
trivial in all degrees and $u_\ast:=\rmd{Id}{Z_\ast}-(-1)^\ast
t_\ast$ is a morphism of complexes from $\mathrm{Z}'_\ast$ to
$\mathrm{Z}_\ast$ (see the definition of Tsygan's double complex).

We claim that $u_\ast$ is injective. Indeed, if
$z\in\rmd{Ker}{}(u_\ast)$ then $t_\ast^{\ast+1}(z)=za$ and, on the
other hand $t_\ast^{\ast+1}(z)=z$. Thus $z(1-a)=0$. But the latter
equality is possible if and only if $z=0$ as $\rmd{Z}{\ast}$ is
free and $1-a$ is not a zero-divisor in $A$. Let
$D_\ast:=\rmd{Coker}{}(u_\ast)$. Then
\[
 0
 \longrightarrow{\mathrm{Z}'_\ast}
 \stackrel{u_\ast}{\longrightarrow}\mathrm{Z}_\ast
 \stackrel{pr}{\longrightarrow}D_\ast
 \longrightarrow0
\]
is an exact sequence of complexes. Since the homology of
$\mathrm{Z}'_\ast$ is trivial and $\mathrm{Z}_\ast$ is acyclic it
results that $D_\ast$ is acyclic.

Our aim now is to show that $D_\ast$ is a projective resolution of
$\rmd{H}{0}(\rmd{Z}{\ast})$ over $\overline{A}$. First let us
prove that each $\rmd{D}{n}$ is $\overline{A}$--projective. It is
easy to see that $\frac{1}{n+1}\sum_{i=0}^n (-1)^{in}t_n^i$
induces a morphism of $\overline{A}$--modules that splits the
canonical map
\[
\rmd{Coker}{}(\rmd{Id}{Z_n}-t_n^{n+1})\longrightarrow\rmd{Coker}{}(u_n).
\]
Thus $D_n$ is projective, as
$\rmd{Coker}{}(\rmd{Id}{Z_n}-t_n^{n+1})=\mathrm{Z}_n/\mathrm{Z}_n(1-a)=\mathrm{Z}_n\otimes_A
\overline{A}$ is free as a right $\overline{A}$--module. To prove
that $D_\ast$ is a resolution it remains to show that
$\rmd{Coker}{}(\overline{d}_1)=\rmd{H}{0}(Z_\ast)$, where
$\overline{d}_\ast$ is the differential of $\rmd{D}{\ast}$. Since
$d'_1$ is surjective we get
$\mathrm{Im}(u_0)\subseteq\mathrm{Im}(d_1)$, hence
\[\rmd{Coker}{}(\overline{d}_1)=[\mathrm{Im}(d_1)+\mathrm{Im}(u_0)]/
\mathrm{Im}(u_0)=\rmd{H}{0}(\mathrm{Z}_\ast).
\]
As $\rmd{char}{}k=0$ we can use Connes' complex
\[
\rmd{C}{\ast}^\lambda:=
\Coker\left(u_\ast\otimes_A{M}\right)=\Coker(u_\ast)\otimes_A{M}
\]
to compute the cyclic homology of $\mathrm{Z}_{\ast}\otimes_A
{M}$. Thus $\rmd{HC}{\ast}(\mathrm{Z}_\ast\otimes_A{M})$ is the
homology of the complex $D_\ast\otimes_A{M}$, where the
differentials of $D_\ast:=\Coker(u_\ast)$ are induced by $d_\ast$.
Finally we have
 \[\rmd{HC}{\ast}(\mathrm{Z}_\ast\otimes_A {M})
 =\rmd{H}{\ast}(\mathrm{D}_\ast\otimes_AM)
 =\rmd{H}{\ast}(\mathrm{D}_\ast\otimes_{\overline{A}}M)
 =\Tor{\overline{A}}{\ast}(\rmd{H}{0}(\mathrm{Z}_\ast),{M}),
 \]
therefore the lemma is proved.
\end{proof}

\begin{noname}\label{Decomp}
Let $K$ be a pointed Hopf algebra, and let $M$ be a right
$K$--comodule. By definition, the coefficient space of $M$ is the
smallest subcoalgebra $C(M)$ of $K$ such that $\rho(M)\subseteq
M\otimes C(M)$. If $C(M)$ is cosemisimple, i.e. it is a (direct)
sum of simple subcoalgebras, then $C(M)=\bigoplus_{x\in{X}}kx$, as
$K$ is pointed. Moreover, in this case, we have
\[\textstyle
M=\bigoplus_{x\in X}M_x,
\]
where $M_x$ contains all $m\in M$ such that $\rho(m)=m\otimes x$.
\end{noname}

\begin{theorem}\label{te:General}
Let $H$ be a Hopf algebra over a field of characteristic $0$ and
let $K$ be a pointed cocommutative Hopf subalgebra of $H$ such
that $H$ is a flat right $K$--module. Let $M\in\mathcal{CM}_m(K)$
such that $C(M)$ is cosemisimple and
$C(M)\subseteq\mathcal{Z}(K)$. Then
\[\textstyle
 \rmd{HC}{\ast}(H,\mathrm{Ind}^H_KM)=
\left(
\bigoplus_{x\in{X_f}}\bigoplus_{i\geq0}\Tor{\overline{K}_x}{\ast-2i}(k,M_x)\right)
 \bigoplus
\left(
\bigoplus_{x\in{X\setminus{X_f}}}\Tor{\overline{K}_x}{\ast}(k,M_x)\right),
\]
where $X_f=\{x\in X\mid \mathrm{ord}\,x<\infty\}$,
$\overline{K}_x=K/(1-x)K$ and $M=\bigoplus_{x\in  X} M_x$ is the
decomposition of $M$ as in \emph{(\ref{Decomp})}.
\end{theorem}

\begin{proof}
Let $x\in X$. Then $\rho(m)=m\otimes{x}$, for every $m\in{M_x}$.
Since $K$ is cocommutative and $C(M)\subseteq{\mathcal{Z}(K)}$, it
results that $M_x$ is a modular crossed $K$--module. By
Proposition~\ref{pr:229} we have
\[\textstyle
 \rmd{HC}{\ast}(H,\mathrm{Ind}^H_KM)\simeq\rmd{HC}{\ast}(K,M)=\bigoplus_{x\in{X}}\rmd{HC}{\ast}(K,M_x).
\]
If $x\in{X_f}$, by Corollary~\ref{co:232}, it follows that
$\rmd{HC}{\ast}(K,M_x)\simeq\bigoplus_{i\geq0}\Tor{\overline{K}_x}{\ast-2i}(k,M_x)$.
In the case when $x\in{X\setminus{X_f}}$, we can apply
Proposition~\ref{pr:030407} to obtain the required description of
$\rmd{HC}{\ast}(K,M_x)$.
\end{proof}

\begin{corollary}\label{co:pointed}
Let $H$ be a pointed Hopf algebra over a field of characteristic
$0$. Let $K:=kG$ be the coradical of $H$ and let
$M\in\mathcal{CM}_m(K)$.

a) There is a set $X\in G$ and a decomposition $M=\bigoplus_{x\in
X}M_x$ as in \emph{(\ref{Decomp})}.

b) If $X$ is as above and, in addition, $X$ is contained in the center of $G$ then
\[\textstyle
 \rmd{HC}{\ast}(H,\rmd{Ind}{K}^HM)=
\left(
\bigoplus_{x\in{X_f}}\bigoplus_{i\geq0}\Tor{\overline{K}_x}{\ast-2i}(k,M_x)\right)
 \bigoplus
\left(
\bigoplus_{x\in{X\setminus{X_f}}}\Tor{\overline{K}_x}{\ast}(k,M_x)\right),
\]
where $X_f=\{x\in X\mid \mathrm{ord}\,x<\infty\}$,
$\overline{K}_x=K/(1-x)K$.
\end{corollary}

\begin{proof}
a) Since $K$ is cosemisimple every right $K$--comodule is a direct
sum of simple comodules, which are one--dimensional as $H$ is
pointed. Furthermore, for every one--dimensional comodule $N$
there is a group-like element $x\in G$ such that $N\simeq kx$.
Thus we can take $M_x$ to be the isotypical component of $M$
corresponding to the simple $K$--comodule $kx$.

b) This part is a direct consequence of the previous theorem.
\end{proof}

Recall that, for an arbitrary group $G$, we denote the set of
conjugacy classes in $G$ by $T(G)$. Let $t(g)$ be a coset of
$T(G)$ and let $t_f(G):=\{x\in t(G)\mid \mathrm{ord}(x)<\infty\}$.
The group homology of $G$ with coefficients in $A$ will be denoted
by $\mathrm{H}_ast(G,A)$.

\begin{corollary}\label{co:group-algebras}
Let $H:=kG$ be the group algebra of an arbitrary group $G$ over a
field $k$ of characteristic $0$. Let $M\in\mathcal{CM}_m(H)$. Then
\[\textstyle
 \rmd{HC}{\ast}(H,M)\simeq
 \left(\bigoplus_{x\in{t_f(G)}}\bigoplus_{i\geq0}\rmd{H}{\ast-2i}(\overline{G}_x,M_x)\right)
 \bigoplus
 \left(\bigoplus_{x\in{t(G)\setminus t_f(G)}}\rmd{H}{\ast}(\overline{G}_x,M_x)\right),
\]
where $M=\bigoplus_{x\in G}M_x$ is the decomposition of $M$ as in
Example \emph{\ref{ex:mcm}(f)}, $G_x$ is the centralizer of $x$ in
$G$, i.e. $G_x:=\{g\in G\mid gx=xg\}$, and
$\overline{G}_x=G_x/\langle x \rangle$.
\end{corollary}
\begin{proof}
By Proposition \ref{pr:Crossed_G} we have
\[\textstyle
 M=\bigoplus_{x\in{t(G)}}\mathrm{Ind}^{kG}_{kG_x}M_x.
\]
Since cyclic homology commutes with direct sums of crossed modules we get
\[\textstyle
\rmd{HC}{\ast}(H,M)\simeq\bigoplus_{x\in{t(G)}}\rmd{HC}{\ast}(kG,\mathrm{Ind}^{kG}_{kG_x}M_x).
\]
To end the proof we apply Corollary \ref{co:pointed}(b).
\end{proof}
\begin{definition}
Let $G$ be a group and let $A$ be a $k$--algebra. $A$ is called
$G$--strongly graded if is a direct sum of $k$--subspaces
$A=\bigoplus_{x\in G}A_x$ such that $A_xA_y=A_{xy}$, $\forall
x,y\in G$.
\end{definition}
\begin{corollary}\label{co:graded-case}
Let $A=\bigoplus_{x\in G}A_x$ be a $G$--strongly graded
$k$--algebra. If $B:=A_1$ and $\overline{A}_x=A_x/[A_x,B]$ then
\[\textstyle
 \rmd{HC}{\ast}(A/B)\simeq
 \left(\bigoplus_{x\in{t_f(G)}}\bigoplus_{i\geq0}\rmd{H}{\ast-2i}(\overline{G}_x,\overline{A}_x)\right)
 \bigoplus
 \left(\bigoplus_{x\in{t(G)\setminus t_f(G)}}\rmd{H}{\ast}(\overline{G}_x,\overline{A}_x)\right).
\]
\end{corollary}

\begin{proof}
We know that $\rmd{HC}{\ast}(A/B)\simeq\rmd{HC}{\ast}(H,A_B)$, see Theorem
\ref{te:S_Bullet}. Since we have $A_B=\bigoplus_{x\in G}\overline{A}_x$ the required
isomorphism follows by Corollary (\ref{co:group-algebras}).
\end{proof}
\begin{remark}
a) If $A=\bigoplus_{x\in G}A_x$ is a $G$--strongly graded algebra
such that $B$ is separable then $\rmd{HC}{\ast}(A/B)$ is
isomorphic to $\rmd{HC}{\ast}(A)$, the usual cyclic homology of
$A$ (of the extension $A/k$). For example, if $B=k$, then we
obtain
\begin{equation}\textstyle\label{re:strong}
 \rmd{HC}{\ast}(A)\simeq
 \left(\bigoplus_{x\in{t_f(G)}}\bigoplus_{i\geq0}\rmd{H}{\ast-2i}(\overline{G}_x,{A}_x)\right)
 \bigoplus
 \left(\bigoplus_{x\in{t(G)\setminus t_f(G)}}\rmd{H}{\ast}(\overline{G}_x,{A}_x)\right).
\end{equation}

b) Let $k$ be a field of characteristic $0$ and let $G$ be an
arbitrary group. The cyclic homology of $kG$ was computed by
Burghelea in \cite{Bu}. By taking $A:=kG$ in (\ref{re:strong}) we
obtain the same result, as $kG$ is $G$--strongly graded with
respect to the decomposition $kG=\bigoplus_{x\in G}kx$.

c) Obviously, if $G$ is a group and $H\unlhd G$ is a normal
subgroup in $G$ then $A:=kG$ is $G/H$--strongly graded, with
$B:=kH$.  In this setting Corollary \ref{co:graded-case} gives the
computation of cyclic homology of the extension $kG/kH$. For a
different approach to this calculation see \cite{Scha}.
\end{remark}

\begin{noname}
An important class of $G$--strongly graded algebras is the class
of $G$--crossed products. Recall that $A=\bigoplus_{g\in G}A_g$ is
called a crossed product if and only if for every $x\in G$ there
is an invertible element $e_x$ in $A_x$. Let $B:=A_1$ and let
$\overline{A}_x=A_x/[A_x,B]$. It is well-known \cite{NvO} that
each $A_x=Be_x$, so $\{e_x\mid x\in G\}$ is a basis of $A$ as a
left $B$ module. We can assume that $e_1$ is the unit of $A$, and
hence of $B$. Define:
\begin{eqnarray}
&\omega:G\times G\to U(B),\qquad\omega(x,y)=e_xe_ye_{xy}^{-1},\label{2cocilu}\\
&``\;.\;":G\times B\to B,\qquad g.b=e_xbe_x^{-1},\label{weak}\\
&\lambda:G\times G\to U(B), \qquad
\lambda(x,y)=e_xe_ye_x^{-1}e_{xyx^{-1}}^{-1}.
\end{eqnarray}
We shall say that $\omega$ is the noncommutative $2$--cocycle of
$A$ and that $``."$ defines the weak action of $G$ on $B$. Then
the multiplication in $A$ is uniquely determined by
\begin{equation}
e_xe_y=\omega(x,y)e_{xy},\qquad e_xb=(g.b)e_x,
\end{equation}
and the action $G_x\times \overline{A}_x\to \overline{A}_x$,
induced by the Ulbrich--Miyashita action of $G$ on $A_B$, is given
by
\[
 g.(\overline{be_x})=\overline{(g.b)\lambda(g,x)e_x}.
\]
Furthermore, we can identify $\overline{A}_x$ with
$\overline{B}_x:=B/[B,B]_x$, where $[B,B]_x$ is the subspace
spanned by all $x$--commutators $bb'-b'(x.b)$.

There are two extreme cases. The first one corresponds to a
trivial cocycle $\omega$, i.e. $\omega(x,y)=1$, for all $x,y\in
G$. One can check easily that (\ref{weak}) defines a real action
of $G$ on $B$ by algebra automorphisms, and
\begin{equation}\label{smash}
(be_x)\cdot(b'e_y)=b(x.b')e_{xy}.
\end{equation}
Therefore $A$ is the smash product $B\#kG$ of $B$ by $G$.

Note that through the identification
$\overline{A}_x\simeq\overline{B}_x$ the action of $G_x$ on
$\overline{A}_x$ corresponds to the action of $G_x$ on $B/[B,B]_x$
induced by $G\times B\to B$.

Conversely, let us assume that $B$ is a $k$--algebra and that $G$
is a group that acts on $B$ by algebra automorphism.  Then
$B\#kG$, the left $B$--module freely generated by $\{e_x\mid x\in
G\}$, becomes an unitary associative algebra with the
multiplication defined in (\ref{smash}). In \cite{FT,GJ}, these
algebras are simply called  crossed products.

In conclusion we have proved the following result.
\end{noname}

\begin{corollary}
Let $B$ be an algebra and let $G$ be a group that acts by
automorphism on $B$. Then:
\[\textstyle
 \rmd{HC}{\ast}(B\#kG/B)\simeq
 \left(\bigoplus_{x\in{t_f(G)}}\bigoplus_{i\geq0}\rmd{H}{\ast-2i}(\overline{G}_x,\overline{B}_x)\right)
 \bigoplus
 \left(\bigoplus_{x\in{t(G)\setminus t_f(G)}}\rmd{H}{\ast}(\overline{G}_x,\overline{B}_x)\right).
\]
If $B$ is a separable algebra then the above isomorphism holds for
$\rmd{HC}{\ast}(B\#kG)$ too.
\end{corollary}

\begin{noname}\label{nn:torus}
Let as consider now the other case, when the weak action
(\ref{weak}) is trivial. Since we are interested in the
computation of the  cyclic homology of quantum tori we shall also
assume that $G=\mathbb{Z}^r$ and $B=k$. Under this assumptions
$\omega$ is a real 2--cocycle of $\mathbb{Z}^r$ with coefficients
in $k^\times$, the multiplicative group of non--zero elements in
$k$, on which $\mathbb{Z}^r$ acts trivially. Moreover
$\overline{A}_x=k$, since every $x$--commutator is zero, and $g\in
\mathbb{Z}^r$ acts on $e_x \in k e_x$ by
\[
 g.e_x={\lambda(g,x)e_x}.
\]
It is easy to see that $g\longmapsto \lambda(g,x)$ defines a
morphisms of groups $\lambda_x:\mathbb{Z}^r\to k^\times$, for
every $x\in\mathbb{Z}^r$.
\end{noname}
\begin{lemma}\label{le:homologyGroups}
Let $G=G_1\times{G_2}$ be a group and let
$\lambda:G\longrightarrow{k^\times}$ be a morphism of groups. Let
$\lambda_i$ denote the restriction of $\lambda$ to $G_i$, $i=1,2$.
If $_\lambda{k}$ is the one--dimensional representation induced by
$\lambda$, then:
\[
 H_n(G,{}_\lambda{k})\simeq\oplus_{p+q=n}H_p(G_1,{}_{\lambda_1}k)\otimes_kH_q(G_2,{}_{\lambda_2}k).
\]
\end{lemma}
\begin{proof}
Obviously we have an isomorphism of $k[G_1\times{G_2}]$--modules:
\[
 _\lambda{k}\simeq{}_{\lambda_1}k\otimes_k{}_{\lambda_2}k,
\]
where the action of $G_1\times{G_2}$ on $_{\lambda_1}k\otimes_k{}_{\lambda_2}k$ is
``{componentwise}". Let $P_\ast\longrightarrow{}_{\lambda_1}k$ and
$Q_\ast\longrightarrow{}_{\lambda_2}k$ be projective resolutions of $_{\lambda_1}k$ and
$_{\lambda_2}k$ over $kG_1$ and $kG_2$ respectively. Then
\[
 \Tot(P_\ast\otimes{Q_\ast})\longrightarrow{}_{\lambda_1}k\otimes_k{}_{\lambda_2}k\simeq{}_{\lambda}k
\]
is a projective resolution of $_\lambda{k}$ over $k[G_1\times{G_2}]$. Proceeding as in
the proof of \cite[Proposition 6.1.13]{We} we get the required isomorphism from the
K{\"u}nneth formula for complexes.
\end{proof}

\begin{corollary}\label{pr:omolog_G}
If $G$ is a finitely generated abelian group and
$\lambda:G\longrightarrow{k^\times}$ is a non--trivial character
of $G$, then $H_\ast(G,{}_\lambda{k})=0$.
\end{corollary}
\begin{proof}
By the preceding Lemma it is enough to show that
$H_\ast(\mathbb{Z},{}_\lambda{k})=0$ and
$H_\ast(\mathbb{Z}_d,{}_\lambda{k})=0$ for every non--trivial
character of $\mathbb{Z}$ and $\mathbb{Z}_d$ respectively. The
equalities follow easily by the computation of group homology for
free and cyclic groups, see \cite[Chapter 6.2]{We}.
\end{proof}

\begin{noname}
Let us keep the notation and assumptions from (\ref{nn:torus}). We
define $A:=k_\omega[\mathbb{Z}^r]$ as follows: as a vector space
$A$ has a basis $\{e_g\mid\;g\in{\mathbb{Z}^r}\}$ and the
multiplication of two elements in the basis is given by
\[
 e_x\cdot{e_y}=\omega(x,y)e_{x+y}.
\]
Obviously $A$ is $\mathbb{Z}^r$--strongly graded, and the
2--cocyle associated to $A$ as in (\ref{2cocilu}) is exactly
$\omega$. Thus for the computation of $HH_{\ast}(A)$ and
$HC_{\ast}(A)$ we can apply Corollary \ref{co:graded-case}. It
results
\[\textstyle
 HC_\ast(A)\simeq
 \bigoplus_{i\geq
 0}H_{\ast-2i}({\mathbb{Z}}^r,ke_0)\bigoplus\left(
 \bigoplus_{x{\in}{\mathbb{Z}^r}{\setminus\{0\}}}H_\ast(\mathbb{Z}^r/\langle x\rangle,ke_x)\right),
\]
where the action of $\mathbb{Z}^r/\langle x\rangle$ on $ke_x$ is
induced by $\lambda_x$. We introduce the notation:
\[
X_\lambda=\{x\in{\mathbb{Z}^r}\mid\lambda_x\mbox{ is
trivial}\}\qquad X_\lambda^\ast= X_\lambda\setminus\{0\}
\]
By Corollary \ref{pr:omolog_G} we get:
\[\textstyle
 HC_\ast(A)=\bigoplus_{i\geq
 0}H_{\ast-2i}({\mathbb{Z}}^r,k)\bigoplus
 \left(\bigoplus_{x\in{X_\lambda^\ast}}H_\ast(\mathbb{Z}^r/\langle
 x\rangle,k)\right).
\]
For each $x=(x_1,\dots,x_n)\in X_\lambda^\ast$ let
$d_x:=\mathrm{gcd}(x_1,\dots,x_n)$. Then
\[\textstyle
\mathbb{Z}^r/\langle
x\rangle\simeq\mathbb{Z}^{r-1}\bigoplus\mathbb{Z}_{d_x}.
\]
Since $k$ is a field of characteristic zero and, for every
$x\in{X_\lambda^\ast}$, the action of $\mathbb{Z}^r/\langle
x\rangle$ on $k$ is trivial, we get
\[\textstyle
 H_\ast(\mathbb{Z}^r/\langle x\rangle,k)\simeq H_\ast(\mathbb{Z}^{r-1}\bigoplus\mathbb{Z}_{d_x},k)
 \simeq H_\ast(\mathbb{Z}^{r-1},k).
\]
By Lemma \ref{le:homologyGroups} it follows that, for every $m\geq
1$, $H_\ast(\mathbb{Z}^m,k) \simeq\Lambda^\ast{V_m}$ where $V_m$
is a vector space of dimension $m$. Hence
\[\textstyle
 HC_\ast(A)\simeq\bigoplus_{i\geq 0}\Lambda^{\ast-2i}{V_r}\bigoplus \left(\Lambda^{\ast}{V_{r-1}}\right)^{(X_\lambda^\ast)}.
\]
Similarly, from the decomposition
$HH_\ast(A)=\bigoplus_{x\in{\mathbb{Z}^r}}H_\ast(\mathbb{Z}^r,ke_x)$,
we deduce:
\[
 HH_\ast(A)=\left(\Lambda^\ast{V_r}\right)^{(X_\lambda)}.
\]
\end{noname}
\noindent Hence we have the following theorem.
\begin{theorem}\label{te:Torus}
If
$\omega:\mathbb{Z}^r\times\mathbb{Z}^r\longrightarrow{k^\times}$
is a normalized $2$--cocycle, then
\[\textstyle
 HH_\ast(k_\omega[\mathbb{Z}^r])=\left(\Lambda^\ast{V_r}\right)^{(X_\lambda)}
 \qquad
 HC_\ast(k_\omega[\mathbb{Z}^r])\simeq\bigoplus_{i\geq 0}\Lambda^{\ast-2i}{V_r}\bigoplus
 \left(\Lambda^{\ast}{V_{r-1}}\right)^{(X_\lambda^\ast)}.
\]
\end{theorem}

\begin{remark}
The algebra $k_\omega[\mathbb{Z}^r]$ is a multiparameter torus. If
$\lambda_{ij}:={\lambda}(f_i,f_j)$, where $\{f_1,\ldots,f_r\}$ is
the canonical basis on $\mathbb{Z}^r$, then
$k_\omega[\mathbb{Z}^r]$ is isomorphic to the algebra generated by
$x_1$, \ldots, $x_r$, $x_1^{-1}$, \ldots, $x_r^{-1}$ satisfying
the following relations:
\[
 x_ix_j=\lambda_{ij}x_jx_i,
 \qquad
 x_ix_i^{-1}=x_i^{-1}x_i=1,
 \qquad
 i, j\in\{1,\ldots,r\}
\]
Hochschild homology of quantum tori is also computed in
\cite{GG,Wa}.
\end{remark}

\section{Cyclic homology of enveloping algebras}
\markboth{\sc Cyclic homology of enveloping algebras}{\sc Cyclic
homology of enveloping algebras}

Let $\Lie{g}$ be a Lie algebra and let $\Ulie{g}$ be its
enveloping algebra. In this section, for every  modular crossed
module $M$ over $\Ulie{g}$, we shall construct a spectral sequence
converging to $\rmd{HC}{\ast}(\Ulie{g},M)$. First we construct
such a spectral sequence for an arbitrary Hopf algebra $H$ and an
arbitrary filtered modular crossed module $M$. Then we shall show
that every crossed module over $\Ulie{g}$ has a canonical
filtration.

Finally, we show that the cyclic homology of almost symmetric
algebras can be computed as the the cyclic homology of $\Ulie g$
with coefficients in a certain modular crossed $\Ulie g$--module,
where $\Lie g$ is a suitable Lie algebra.
\begin{theorem}
Let $M$ be a filtered modular crossed module with increasing
filtration $(F_pM)_{p\in\mathbb{N}}$. Then there is a spectral
sequence
\begin{equation}\label{eq:spectral}
\rmd{E}{pq}^1=\rmd{HC}{p+q}(gr_pM)\Longrightarrow\rmd{HC}{p+q}(M),
\end{equation}
where $gr_p(M)=F_pM/F_{p-1}M$, for every $p\geq0$, ($F_{-1}M=0$).
\end{theorem}
\begin{proof}
The spectral sequence that we construct is a particular case of
the spectral sequence associated to a filtration.

By definition $\rmd{HC}{\ast}(H,M)$ is the homology of
$\rmd{Tot}{}(\rmd{CC}{\ast\ast}(H,M))$. Since the functor from
$\mathcal{CM}_m(H)$ to the category of chain complexes
\[
 M\longmapsto\rmd{Tot}{}\rmd{CC}{\ast\ast}(H,M)
\]
is exact, we have an exact sequence
\[
 0\longrightarrow
 \rmd{Tot}{}\rmd{CC}{\ast\ast}(H,F_{p-1}M)\rightarrow
 \rmd{Tot}{}\rmd{CC}{\ast\ast}(H,F_{p}M)\rightarrow
 \rmd{Tot}{}\rmd{CC}{\ast\ast}\left(H,\frac{F_{p}M}{F_{p-1}M}\right)\rightarrow0.
\]
As $\rmd{CC}{\ast\ast}(H,-)$ commutes with direct sums, the graded
complex associated to the filtration
$\left(\rmd{Tot}{}\rmd{CC}{\ast\ast}(H,F_pM)\right)_{p\geq0}$ is
$\bigoplus_{p\geq0}\rmd{Tot}{}\rmd{CC}{\ast\ast}\left(H,\frac{F_pM}{F_{p-1}M}\right)$.
It follows that the spectral sequence associated to this
filtration has the following $0$--page
\[
 E^0_{pq}
 =
 \frac{F_p\rmd{Tot}{p+q}\rmd{CC}{\ast\ast}(H,M)}
      {F_{p-1}\rmd{Tot}{p+q}\rmd{CC}{\ast\ast}(H,M)}
 =
 \rmd{Tot}{p+q}\rmd{CC}{\ast\ast}\left(H,\frac{F_pM}{F_{p-1}M}\right).
\]
We conclude by remarking that
$\rmd{HC}{\ast}\left(H,\frac{F_pM}{F_{p-1}M}\right)
=\rmd{H}{\ast}\left(\rmd{Tot}{}\rmd{CC}{\ast\ast}\left(H,\frac{F_pM}{F_{p-1}M}\right)\right)$.
\end{proof}

\begin{lemma}
For each crossed module $M$ over $\Ulie{g}$ there is a
non--negative increasing filtration $(F_pM)_{p\in\mathbb{Z}}$ on
$M$ such that $\bigcup_{p\in\mathbb{Z}}F_pM=M$ and $F_pM/F_{p-1}M$
is a comodule with trivial structure for every $p\in\mathbb{Z}$.
\end{lemma}

\begin{proof}
Recall that $M^{coH}=\{m\in M\mid \rho(m)=m\otimes 1\}$. We set
$F_pM=0$, if $p<0$, and $F_{0}M=M^{coH}$. Inductively we define
$F_{p+1}M$ such that
\[
F_{p+1}M/F_{p}M=(M/F_{p}M)^{coH}.
\]
Obviously the comodule structure of the graded associated to this
filtration is trivial. Since $\Ulie{g}$ is cocommutative it
results easily that each $F_{p}M$ is a submodule of $M$. It
remains to prove that $\bigcup_{p\in\mathbb{N}}F_pM=M$. Let
\[\textstyle
N:=\bigcup_{p\in\mathbb{N}}F_pM
\]
and let us assume that $N\neq M$. Then $M/N$ contains a simple
$\Ulie{g}$--subcomodule. Since there is only one type of simple
$\Ulie{g}$--comodules and that one has trivial comodule structure
($\Ulie{g}$ is a connected Hopf algebra, i.e. its coradical is one
dimensional) it follows that $(M/N)^{coH}\neq 0$. Let
$\overline{m}\in (M/N)^{coH}$ be a non--zero element. Thus:
\[
\rho(m)=m\otimes 1+x,
\]
where $x\in N\otimes H$. As $x$ is a finite sum of tensor
monomials and $N=\bigcup_{p\in\mathbb{N}}F_pM$ it follows that
there is $p_0$ such that $x\in F_{p_0}M\otimes H$. Hence
$\rho(m)=m\otimes 1+F_{p_0}M\otimes H$. Thus $m\in F_{p_0+1}M$,
contradicting the fact that the class $\overline{m}\in M/N$ is not
zero.
\end{proof}

\begin{theorem}
Let $M$ be modular crossed module over $\Ulie{g}$. Let
$(F_pM)_{p\in\mathbb{N}}$ be the filtration constructed in the
previous lemma. Then there is a spectral sequence
\[\textstyle
\rmd{E}{pq}^1=\bigoplus_{i\geq0}\rmd{H}{p+q-2i}(\Lie{g},F_pM/F_{p-1}M)\Longrightarrow
\rmd{HC}{\ast}(M).
\]
\end{theorem}

\begin{proof}
By Theorem \ref{th:228} we have
\[\textstyle
\rmd{HC}{p+q}(\Ulie{g},F_pM/F_{p-1}M)\simeq\bigoplus_{i\geq0}\rmd{H}{p+q-2i}(\Lie{g},F_pM/F_{p-1}M).
\]
Hence the spectral sequence that we are looking for is a
particular case of (\ref{eq:spectral}).
\end{proof}

We shall end the paper by studying the cyclic homology of a
certain class of filtered algebras, that appeared in the work of
C. Kassel.
\begin{definition}
Let $A$ be a filtered algebra with non--negative increasing
filtration $(F_nA)_{n\in\mathbb{N}}$. Following \cite[p. 100]{Lo}
we shall say that $A$ is an almost symmetric algebra if $\gr{A}$,
the graded associated to $A$, is isomorphic to $S(V)$, the
symmetric algebra of a vector space $V$.
\end{definition}
Sridharan classified all almost symmetric algebras in \cite{Sr}.
He proved that for each almost symmetric algebra $A$ there are a
Lie algebra $(\Lie{g},[-,-])$ and a 2--cocycle
$f:\Lie{g}\otimes\Lie{g}\to k$ such that $\Lie g\simeq F_1A/F_0A$
(as vector spaces) and
\[
A\simeq T(\Lie{g})/I_f(\Lie g),
\]
where $I_f(\Lie g)=\langle x\otimes y-y\otimes x- [x,y]-f(x\otimes y)\mid x,y\in
\Lie{g}\rangle$. If we consider the standard filtration on $T(\Lie g)$ then the above
isomorphism respects the filtrations.

Let  $\Uf g:=T(\Lie g)/I_f(\Lie g)$. There is a canonical
$k$--linear map $i_f:\Lie g\to\Uf g$. For every $x\in \Lie g$ the
element $i_f(x)$ will be denoted by $\overline x$. The proof of
the above description of $A$ is based on the fact that $\Uf g$ has
a PBW-basis. Suppose that $(x_i)_{i\in I}$ is a basis of $\Lie g$,
with $I$ a totally ordered set. Then the set containing $1$ and
all monomials $\overline x_{i_1}\overline x_{i_2}\cdots \overline
x_{i_n}$, with $i_1\leq i_2\leq\ldots\leq i_n$, is a basis on $\Uf
g$.

Furthermore, one can construct an $\Ulie{g}$--comodule algebra
structure on $\Uf g$ in the following way. First we define an
algebra map $T(\Lie{g})\longrightarrow A\otimes \Ulie{g}$ such
that
\[
x\longmapsto \overline{x}\otimes 1+1\otimes x, \forall x\in
\Lie{g}.
\]
It factorizes through an algebra map $\rho:A\to A\otimes\Ulie{g}$.
One easily checks that $(A,\rho)$ is an $\Ulie{g}$--comodule
algebra such that
\begin{equation}\label{UfCoaction}
\overline{x}\longmapsto \overline{x}\otimes 1+1\otimes x, \forall
x\in \Lie{g}.
\end{equation}

Moreover, the set whose elements are $1$ and all monomials $x_{i_1}x_{i_2}\cdots
x_{i_n}$, with $i_1\leq i_2\leq\ldots\leq i_n$, is a basis on $\Ulie g$. Thus there is
a unique $k$--linear map $\theta:\Ulie g\to\Uf g$ such that $\theta(1)=1$ and, for
every $i_1\leq i_2\leq\ldots\leq i_n$, we have
\[
\qquad\theta(x_{i_1}x_{i_2}\cdots x_{i_n})=\overline x_{i_1}\overline x_{i_2}\cdots
\overline x_{i_n}.
\]
A straightforward computation shows us that $\theta$ is $H$--colinear. Recall that, by
definition, a total integral for an $H$--comodule algebra $A$ is an $H$--colinear map
from $H$ to $A$ that maps the unit of $H$ to the unit of $A$. Thus $\theta$ is a total
integral for $\Uf g$.

Since $\theta$ is bijective too it follows that the subalgebra of
coinvariant elements in $\Uf g$ is $k$. Then, by \cite[Proposition
1.5]{Bell}, it results that the extension $\Uf g/k$ is $\Ulie
g$--Galois and $\theta$ is invertible in convolution, i.e. there
is a $k$--linear function $\theta^{-1}:\Ulie g\to \Uf g$ such that
\[\textstyle
\sum\theta(h\de 1)\theta^{-1}(h\de 2)=\varepsilon(h)1_A=\theta^{-1}(h\de 1)\theta(h\de
2),\qquad \forall h\in\Ulie g.
\]
Hence, by the above relation, it results that  $\kappa(h)=\tsum \theta^{-1}(h\de
1)\otimes\theta(h\de 2)$, see (\ref{pa:Galois}) for the definition of $\kappa$. Thus
the Ulbrich--Miyashita action of $\Ulie g$ on $\Uf g$ is uniquely determined by
\begin{equation}\label{Ulbr}
x.a=\overline{x}a-a\overline{x},\qquad\forall x\in\Lie g,\,\,\forall a\in\Uf g.
\end{equation}
Summarizing, we have the following result.
\begin{proposition}
Let $\Lie g$ be a Lie algebra and let $f:\Lie g\otimes\Lie g\to k$
be a $2$--cocycle. Then $\Uf g$ is an $\Ulie g$--Galois extension
of $k$ and
\begin{equation*}\label{UfAction}
\rmd{HC}{\ast}(\Uf g)=\rmd{HC}{\ast}(\Ulie g,\Uf g),
\end{equation*}
where $\Uf g$ is a  modular crossed module with respect to the
action \emph{(\ref{Ulbr})} and coaction \emph{(\ref{UfCoaction})}.
\end{proposition}


\end{document}